\documentclass[smallextended]{svjour3}
\smartqed  
\usepackage{fix-cm}
\usepackage{amssymb}
\usepackage{epsfig}
\usepackage{rotating}
\usepackage{caption}
\usepackage{lscape}
\usepackage{graphics}
\usepackage{xcolor}
\usepackage{algorithm}
\usepackage{algorithmicx}
\usepackage{algpseudocode}
\usepackage{amsmath}
\usepackage{amsfonts}
\usepackage{epic,eepic}
\usepackage{epsfig}
\usepackage[toc,page,title,titletoc,header]{appendix}

\newtheorem{thm}{Theorem}
\newtheorem{lem}{Lemma}
\newtheorem{rem}{Remark}
\newtheorem{asmp}{Assumption}
\newtheorem{exam}{Example}
\newtheorem{dfn}{Definition}
\renewcommand\endproof{\hfill $\Box$ \vskip .5cm}

\begin{document}

\title{SOCP Reformulation for the Generalized Trust Region Subproblem via a Canonical Form of Two Symmetric Matrices}

\titlerunning{SOCP Reformulation for GTRS}        

\author{Rujun Jiang \and Duan Li \and Baiyi Wu }


\institute{
          Rujun Jiang \at Department of Systems Engineering and
Engineering Management, The Chinese University of Hong Kong, Shatin,
N. T., Hong Kong \\
              \email{rjjiang@se.cuhk.edu.hk}
               \and  Duan Li \at
              Corresponding author. Tel.: +852-39438323, Fax: +852-26035505\\Department of Systems Engineering and
Engineering Management, The Chinese University of Hong Kong, Shatin,
N. T., Hong Kong \\
              \email{dli@se.cuhk.edu.hk}
              \and  Baiyi Wu \at
              Department of Systems Engineering and
Engineering Management, The Chinese University of Hong Kong, Shatin,
N. T., Hong Kong \\
              \email{baiyiwu@gmail.com}}

\date{Received: date / Accepted: date}

\maketitle

\begin{abstract}
We investigate in this paper the generalized trust region subproblem (GTRS) of minimizing a general  quadratic objective function subject to a general quadratic inequality constraint.
By applying a simultaneous block diagonalization approach,  we obtain a congruent canonical form
for the symmetric matrices in both  the objective and  constraint functions.
By exploiting the block separability of the canonical form, we show that  all GTRSs with an optimal value bounded from below are second order cone programming (SOCP) representable.
Our result generalizes the recent work of  Ben-Tal and Hertog  (Math. Program. 143(1-2):1-29, 2014), which  establishes the SOCP representability of the GTRS under the assumption of the simultaneous diagonalizability of the two matrices in the objective and constraint functions.
 Compared with the state-of-the-art approach to reformulate the GTRS as a semi-definite programming  problem, our SOCP reformulation delivers
 a much faster solution algorithm. We further extend our method to two variants of the GTRS in which the inequality constraint is replaced by either an equality constraint or an interval constraint.
 Our methods also enable us to obtain simplified versions of the classical S-lemma, the S-lemma with equality, and the S-lemma with interval bounds.
\keywords{Trust region subproblem \and Simultaneous block diagonalization \and
Quadratically constrained quadratic programming \and  Second order cone programming}

 \subclass{ 90C20\and  90C26}
\end{abstract}

\section{Introduction}
We consider in this paper the following generalized trust region subproblem (GTRS):
\begin{eqnarray*}
\rm{(P)}~
&\min
&f(x)=\frac{1}{2}x^TDx+e^Tx \\
&{\rm s.t.}& h(x)=\frac{1}{2}x^TAx+b^Tx+c\leq 0,
\end{eqnarray*}
where  $A$ and $D$ are $n\times n$ symmetric matrices but not necessarily positive semi-definite, $b,e\in \Re^n$ and $c\in \Re$.

The GTRS has been widely investigated in the optimization literature and includes the classical trust region subproblem as its special case
where the constraint reduces to a ball constraint, i.e., $x^Tx\le 1$. The classical trust region subproblem is fundamental in trust region methods for nonlinear optimization problems, see
\cite{more1993generalizations,stern1995indefinite,yuan1990subproblem}.
Other prominent applications of the classical trust region subproblem can be found in regularization and robust optimization \cite{ben2014hidden}, etc. Rendl and Wolkowicz
\cite{rendl1997semidefinite} first solve the classical trust region subproblem via a semi-definite programming (SDP) reformulation. The past two decades have witnessed
numerous methods developed for solving  the GTRS under various assumptions, see, for example,
\cite{ben1996hidden,feng2012duality,more1993generalizations,pong2014generalized,sturm2003cones,yuan1990subproblem}.
Sturm and Zhang
further reveal that problem ${\rm(P)}$ admits an  SDP reformulation. Still, the relatively large computational complexity of SDP algorithms prevents them from scaling to
large-scale problems.
Most fast algorithms  \cite{more1993generalizations,pong2014generalized} for the GTRS are developed under a regular condition that there exists
a $\lambda\in\Re$ such that  $D +\lambda A\succ0$, together with  some other mild conditions.
  Recently, Ben-Tal and Hertog \cite{ben2014hidden} show that if the two matrices in both the objective and   constraint functions  are
  simultaneously diagonalizable, the GTRS can be then transformed into an equivalent second order cone programming (SOCP) problem formulation. Simultaneous
  diagonalizability is  actually a more general condition than the regular condition \cite{feng2012duality}. Conditions for simultaneous diagonalizability and corresponding algorithms
   are investigated recently in  \cite{jiang2015simultaneous}. Compared with the SDP representation of the GTRS,
   the SOCP representation delivers a much faster solution algorithm, which is
  critical for solving large-scale problems in practice. This recognition motivates the investigation in our study.

We advance the state-of-the-art methods for solving the GTRS in this paper. More specifically, we prove that all GTRSs with an optimal value bounded from below are SOCP representable. To obtain the SOCP representation,
 we invoke and extend the congruent canonical form in Uhlig \cite{uhlig1976canonical}. We first transform the two matrices
  into their canonical form of block diagonal matrices via congruence, and then make use of the block separability of the matrices in the canonical form.
  In particular,  we derive necessary conditions  from the canonical form for the GTRS to be bounded from below, and then show that the problem can further be transformed to an SOCP reformulation under such necessary conditions.  Moreover, the attainableness of the optimal value is determined by the associated coefficients in the linear terms in both the objective and constraint functions.
Note that the transformation into separable quadratic forms (the canonical form of block diagonal matrices) can be done off-line (independent of the algorithm)
 and the time complexity of the transformation is almost negligible when  compared with that of the SDP reformulation. In fact, our method using simultaneous block diagonalization is a generalization of the simultaneous
  diagonalizability in \cite{ben2014hidden}.

We also extend our approach to
   two variants of problem   ${\rm(P)}$ where the constraint is replaced by either an  equality constraint,
\begin{eqnarray*}
\rm{(EP)}~&~ \min & f(x)=\frac{1}{2}x^TDx+e^Tx \\
&~  {\rm s.t.}&h(x)=\frac{1}{2}x^TAx+b^Tx+c=0,
\end{eqnarray*}
or an interval constraint,
\begin{eqnarray*}
\rm{(IP)}~&~ \min & f(x)=\frac{1}{2}x^TDx+e^Tx \\
&~  {\rm s.t.}&c_1\leq h(x)=\frac{1}{2}x^TAx+b^Tx\leq c_2.
\end{eqnarray*} Mor{\'e} \cite{more1993generalizations} presents  a method for problem ${\rm(EP)}$
  by using the saddle point optimality condition
 under some mild assumptions. Xia et al. \cite{xia2014s} transform the problem ${\rm(EP)}$
  to an SDP reformulation by using the S-lemma with equality under the conditions that $A\neq0$
 and $h(x)$ can take both positive and negative values. One application of problem  ${\rm(EP)}$
 is  time of arrival problem
\cite{hamm2010quadratic}. Stern and Wolkowicz \cite{stern1995indefinite}  propose a method for problem ${\rm(IP)}$ under $b=0$ and the regular condition.
 By assuming  $b=0$ and the simultaneous diagonalizability of $A$ and $D$, Ben-Tal  and  Teboulle \cite{ben1996hidden} derive the hidden convexity of
 problem ${\rm (IP)}$ and thus transform the problem to an SOCP reformulation. Ye and Zhang \cite{ye2003new} further show that problem ${\rm(IP)}$
  admits an SDP reformulation if both the primal and dual Slater conditions are satisfied. Recently, strong duality conditions of  ${\rm (IP)}$ are
  studied in Pong and Wolkowicz \cite{pong2014generalized} and a fast method is provided   under the regular condition. Wang and Xia \cite{wang2014strong} further
  simplify the
 conditions in  \cite{pong2014generalized} and develop the S-lemma with interval bounds to solve ${\rm (IP)}$. Ben-Tal and Hertog \cite{ben2014hidden}
  further show that  ${\rm (IP)}$ can be solved as an SOCP when $A$ and $D$ are simultaneously diagonalizable without the assumption of $b=0$. Note that
   ${\rm (IP)}$ includes the equality
 constrained problem ${\rm (EP)}$ as a special case when setting $c_1=c_2$. On the other hand, we will discuss in latter sections that solution methods for  ${\rm (EP)}$ can also be used to solve  ${\rm (IP)}$. Essentially, we will show that some slightly modified versions of our previous results for problem ${\rm(P)}$ hold true for the equality constrained problem ${\rm(EP)}$
   and the interval constrained problem ${\rm(IP)}$.

 To summarize, we derive necessary conditions for problem ${\rm(P)}$ and its variants with equality constraint or interval constraint to be bounded from below
and further transform the problems to their SOCP reformulations by exploiting the block separability of the canonical form. Besides, we also derive the conditions for the attainableness of the problems. We emphasize that our methods are applicable for
 general situations without specific assumptions such as the regular condition. We would also like to indicate that
 our methods avoid involvement of linear matrix inequalities (LMI) as LMIs are usually hard to handle for large-scale problems.
As by-products of our research, we further obtain simplified versions of the classical S-lemma, the S-lemma with equality and the S-lemma with interval bounds.

We organize our paper as follows. In Section 2, we introduce and extend a  canonical form for the two matrices in both the objective and  constraint functions by
a real congruent transformation. After identifying all the situations in which the optimal value is bounded from below, we show that all remaining bounded cases
of problem ${\rm(P)}$ can be transformed into an SOCP reformulation. In Section 3,  we extend our methods to problems ${\rm(EP)}$ and ${\rm(IP)}$. Finally,
we conclude our paper in Section 4.

{\textbf{Notations}}: Throughout  this paper, $I_{m}$ represents the $m\times m$ identical matrix. $\bar{1}$ denotes the all one vector $(1,\ldots,1)^T$. The notation $\Re^n$ represents the $n$ dimensional vector space. For symmetric matrices $A$ and $B$, $A$
$\succeq$ $B$ denotes that matrix $A-B$ is positive semi-definite. We denote the \emph{Moore--Penrose pseudoinverse} by $A^+$. We use $\text{sign}(x)$ to denote the \text{sign} of a real number $x$, i.e., $\text{sign}(x)=1$, if $x\geq0$, otherwise $\text{sign}(x)=-1$. And we use $\dim A$ to denote the dimension of a square matrix $A$. And $A_{k:l,k:l}$ denotes the submatrix of matrix A by selecting the rows $k,k+1,\ldots,l$, and the columns $k,k+1,\ldots,l$. We also denote by  ${\rm diag}(A_1,\ldots,A_k)$
 the block diagonal matrix \[\left(
\begin{array}{ccccc}
A_1&&&0\\
&\cdot&&\\
&&\cdot&&\\
0&&&A_k
\end{array}
\right).\]
We denote by $E$ the anti-diagonal matrix
\begin{equation}\label{refE}
\left(
\begin{array}{ccccc}
0&&&&1\\
&&&\cdot&\\
&&\cdot&&\\
&\cdot&&&\\
1&&&&0
\end{array}
\right),\end{equation}
and by $F$ the lower striped matrix
\begin{equation}\label{refF}\left(\begin{array}{ccccc}
 & & & &0\\
 & & &0&1 \\
&&\ldots&&\\
 &0&1&&\\
0&1 &&&
\end{array}\right).\end{equation}
We use $J(\lambda,m)$ to denote an $m\times m$  Jordan block
\[\left(
\begin{array}{ccccc}
\lambda&e&&&\\
&\cdot&\cdot&&\\
&&\cdot&\cdot&\\
&&&\lambda&e\\
&&&&\lambda
\end{array}
\right).\]
If the eigenvalue is a real number, i.e., $\lambda\in\Re$, then $e=1$ for $m\geq2$, while $J=(\lambda)$ for $m=1$.
If the eigenvalues form a complex pair, i.e., $a \pm bi$, then $\lambda=\left(
\begin{array}{cc}
a&-b\\
b&a
\end{array}
\right)$, with $a$ and $b$ $\in\Re$, and $b\neq0$, and $e=\left(
\begin{array}{cc}
1&0\\
0&1
\end{array}
\right)$ for $m\geq4$, while  $J=\left(
\begin{array}{cc}
a&-b\\
b&a
\end{array}
\right)$ for $m=2$.

Let $v{\rm(P)}$ denote the optimal value of problem ${\rm(P)}$. For an optimization problem $\min\{f(x)\mid x\in X\}$ with a nonempty feasible
set $X$, we define its optimal value as $v=\inf\{f(x)\mid x\in X\}$. For any $\epsilon>0$, we call a solution $\bar{x}\in X$  an $\epsilon$
 optimal solution, if $f(\bar{x})-v\leq\epsilon$. Note that we reuse the  notations $h(x)$ and $f(x)$
in problems  ${\rm(P)}$, ${\rm(EP)}$  and ${\rm(IP)}$ and the simplified versions of the S-lemma with different constraints to keep the notations simple and help
readers follow the main theme of the paper with ease.

\section{ SOCP reformulation for GTRS}
In this section, we first use congruent transformation to transform any two symmetric matrices into a canonical form of two block diagonal matrices.
 Then we explore different cases of the
canonical form with respect to the boundedness of the optimal value and its attainability and transform problem ${\rm(P)}$ with its optimal value bounded from below to an equivalent SOCP problem.
\subsection{Canonical form of two symmetric matrices}
We invoke the following lemmas from Uhlig \cite{uhlig1976canonical} to obtain a canonical form of any two real symmetric matrices.
\begin{lem}\label{L2.1}
(Theorem 1 in \cite{uhlig1976canonical})
Let $A$ and $D$ be two $n\times n$ real symmetric matrices. Suppose $A$ is nonsingular. Let $A^{-1}D$ have a Jordan normal form $diag (J_1,\ldots,J_k)$, where
$J_1,\ldots,J_k$ are Jordan blocks either with real eigenvalues or with complex eigenvalues. Then there exists an $n\times n$ real congruent matrix $S$ such that
\begin{eqnarray*}
S^T AS= {\rm diag}(\varepsilon_1E_1,\ldots,\varepsilon_kE_k)
\end{eqnarray*}
and
\begin{eqnarray*}
S^TDS= {\rm diag}(\varepsilon_1E_1J_1,\ldots,\varepsilon_kE_kJ_k),
\end{eqnarray*}
where $\varepsilon_i=\pm1$ and $E_i$ is defined in (\ref{refE}). Furthermore, the \text{sign}s of $\varepsilon_i$, $i$ = 1, $\ldots$, $k$, are uniquely (up to permutations) associated with the Jordan blocks, $J_i$, $i$ = 1, $\ldots$, $k$. In particular,
 $\varepsilon_i=1$ if $J_i$ has a pair of complex eigenvalues.
\end{lem}

\begin{lem}\label{L2.2}
(Theorem 2 in \cite{uhlig1976canonical}) Let $A$ and $D$ be two singular real symmetric matrices and assume that there exists a $\mu \in\Re$ such  that $C$ = $A+\mu D$
is nonsingular. Let
\begin{eqnarray*}
&J&={\rm diag}(J(\lambda_1,n_1),\ldots,J(\lambda_k,n_k),J(0,n_{k+1}),\ldots,J(0,n_p),J(1/\mu,n_{p+1}),\ldots,\\
&&~~~~~~~~~J(1/\mu,n_m))
\end{eqnarray*}
be the  Jordan normal form of $C^{-1}D$.
Then there exists an $n\times n$ real congruent matrix $S$ such that
\begin{equation}
S^T AS={\rm diag}(\tau_1E_1,\ldots,\tau_k E_k,\tau_{k+1}E_{k+1},\ldots,\tau_pE_p,\tau_{p+1}F_{p+1},\ldots,\tau_mF_m) \label{Equation1}
\end{equation}
and
\begin{eqnarray}
&S^TDS=&{\rm diag}(\tau_1E_1J(\kappa_1,n_1),\ldots,\tau_kE_kJ(\kappa_k,n_k),\tau_{k+1}F_{k+1},\ldots,\tau_pF_p,\notag\\
&&~~~~~\tau_{p+1}E_{p+1},\ldots,\tau_mE_m), \label{Equation2}
\end{eqnarray}
where $E_i$ and $F_i$ are defined in (\ref{refE}) and (\ref{refF}), respectively, $\tau_i=\pm1,~i=1,\ldots,m,$ $\dim{E_i}=\dim{F_i}=n_i$, $i$ = $k+1$, $\ldots$, $m$, and $\kappa_i=\lambda_i/(1-\mu\lambda_i)$, $i$ = 1, $\ldots$, $k$.
The \text{sign}s of $\tau_i$ are uniquely (up to permutations) determined by the associated Jordan blocks $J(\kappa_i,n_i), ~E_i$ or $F_i$. In particular,  $\tau_i=1$ if $J(\lambda_i,n_{i})$ has a pair of complex
eigenvalues, $i=1,\ldots,k$. Furthermore, $p-k\ge1$ and $m-p\ge 1$.
\end{lem}
%
%

Next we generalize the results in the previous two lemmas to general situations where  we do not assume the existence of the nonsingular matrix pencil for two symmetric matrices.
\begin{thm}\label{T2.1}
For any two $n$ $\times$ $n$ real symmetric matrices $A$ and $D$, there exists an $n\times n$ real invertible matrix $S$ such that
\begin{eqnarray}
&S^T AS=&{\rm diag}(\tau_1E_1,\ldots,\tau_kE_k,\tau_{k+1}E_{k+1},\ldots,\tau_pE_p,\tau_{p+1}F_{p+1},\ldots,\tau_mF_m,\notag\\
&&~~~~~~0,\ldots,0)\label{CF1}
\end{eqnarray}
and
\begin{eqnarray}
&S^TDS=&{\rm diag}(\tau_1E_1J(\kappa_1,n_1),\ldots,\tau_kE_kJ(\kappa_k,n_k),\tau_{k+1}F_{k+1},\ldots,\tau_pF_p,\notag\\
&&~~~~~~\tau_{p+1}E_{p+1},\ldots,\tau_mE_m,0\ldots,0)\label{CF2}
\end{eqnarray}
where $\dim{E_i}=\dim{F_i}=n_i,~i=k+1,\ldots,m$, and $\tau_i=\pm1,~i=1,\ldots,m.$ The signs of $\tau_i$ are uniquely (up to permutations) determined by the associated
 Jordan blocks $J(\kappa_i,n_i), ~E_i$ or $F_i$. The values of $\kappa_i$ are uniquely (up to permutations) determined by the associated Jordan blocks $J(\kappa_i,n_i)$.
\end{thm}
\proof
Given Lemmas \ref{L2.1} and \ref{L2.2}, we only need to consider in the the proof  the case where   $A$ and $D$ are both singular and there does not exist a $\mu\in\Re$ such that $A+\mu D$ is nonsingular.

We can always find a congruent matrix $Q_1$ such that  $\bar{A}\triangleq Q_1^TAQ_1={\rm diag}(A_1,0,\ldots,0)$, where $A_1$ is a $q\times q$ diagonal matrix and $q={\rm rank}(A)$. Denote
\[\bar{D}\triangleq Q_1^TDQ_1=
\left(\begin{array}{cc}
D_1&D_2\\
D_2^T&D_3\\\end{array}\right),\]
where $D_1$ is a $q \times q$ matrix. We can always find a congruent matrix $S$ such that $S^TD_3S={\rm diag}(D_6,0,\ldots,0)$, where $D_6$ is a nonsingular $s\times s$ diagonal matrix. Let $Q_2\triangleq {\rm diag}(I_q, S)$.
Then $\hat{A}\triangleq Q_2^T\bar{A}Q_2=\bar{A}$, and
\[\hat{D}\triangleq Q_2^T\bar{D}Q_2=\left(\begin{array}{ccc}
D_1&D_4&D_5\\
D_4^T&D_6&0\\
D_5^T&0&0
\end{array}\right).\]
Let \[Q_3\triangleq\left(\begin{array}{ccc}
I_q&0&0\\
-D_6^{-1}D_4^T&I_{s}&0\\
0&0&I_{n-q-s}
\end{array}\right).\]
Then,
\[\tilde{D}\triangleq Q_3^T\hat{D}Q_3=\left(\begin{array}{ccc}
D_1-D_4D_6^{-1}D_4^T&0&D_5\\
0&D_6&0\\
D_5^T&0&0
\end{array}\right),\] and $\tilde{A}\triangleq Q_3^T\hat{A}Q_3=\hat{A}=\bar{A}$.
We can always choose a $\mu\in\Re$ such that the first $q$ columns of $\tilde{D}+\mu \tilde{A}$ are linearly independent. For example, we can choose $\mu=\max_{i=1,\ldots,q}\sum_{j=1}^q |b_{ij}|/|a_{ii}|+1$, where $b_{ij}$ is the element in the $i$th row and the $j$th column of $D_1-D_4D_6^{-1}D_4^T$ and $a_{ii}$
is the $i$th diagonal element of $A_1$.
Then $\mu A_1+D_1-D_4D_6^{-1}D_4^T$ is nonsingular.

If the columns in $D_5$ are linearly independent, then $\tilde{D}+\mu \tilde{A}$
is nonsingular and thus $D+\mu A$ is nonsingular,  which contradicts our assumption of no nonsingular matrix pencil. Thus the columns in $D_5$ are linearly dependent.
We can always find a congruent matrix $Q_4$ such that
$\check{A}=Q_4^T\tilde{A}Q_4=\hat{A}=\bar{A}$,
and
$$\check{D}\triangleq Q_4^T\tilde{D}Q_4=\left(\begin{array}{cccc}
D_1-D_4D_6^{-1}D_4&0&D_5'&0\\
0&D_6&0&0\\
D_5'^T&0&0&0\\
0&0&0&0
\end{array}\right),$$ where $D_5'$ is of full column rank. Let
$$A'\triangleq\left(\begin{array}{ccc}
A_1&0&0\\
0&0&0\\
0&0&0
\end{array}\right),~
D'\triangleq\left(\begin{array}{ccc}
D_1-D_4D_6^{-1}D_4^T&0&D_5'\\
0&D_6&0\\
D_5'^T&0&0
\end{array}\right).$$
Then there exists a $\mu$ such  that $D'+\mu A'$ is nonsingular.
From Lemma \ref{L2.2} we know that $A'$ and $D'$ can be congruent to the canonical form in (\ref{Equation1}) and (\ref{Equation2}).
So $A$ and $D$ can be congruent to the canonical form in (\ref{CF1}) and (\ref{CF2}).
\endproof

\begin{rem}
If $A$ and $D$ are both singular and there does not exist an $\mu\in\Re$ such that $A+\mu D$ is nonsingular, then the number of  the common $0$ terms in the lower right part of (\ref{CF1}) and (\ref{CF2}) is equal to $n-{\rm rank}(A'+\mu D')$.
\end{rem}

From Lemma \ref{L2.1}, Lemma \ref{L2.2} and Theorem \ref{T2.1}, we know that (\ref{CF1}) and (\ref{CF2}) represent a canonical  form for any two real symmetric matrices $A$ and $D$ via congruence.
Without loss of generality,  we assume from now on that matrices $A$ and $D$ in problem (P) satisfy:
\begin{eqnarray*}
A&=&{\rm diag}(A_1,A_2,\ldots,A_s)\\
~&=&{\rm diag}(\tau_1E_1,\ldots,\tau_kE_k,\tau_{k+1}E_{k+1},\ldots,\tau_pE_p,\tau_{p+1}F_{p+1},\ldots,\tau_mF_m,0,\ldots,0),\\
D&=&{\rm diag}(D_1,D_2,\ldots,D_s)\\
~&=&{\rm diag}(\tau_1E_1J(\kappa_1,n_1),\ldots,\tau_kE_kJ(\kappa_k,n_k),\tau_{k+1}F_{k+1},\ldots,\tau_pF_p,\tau_{p+1}E_{p+1},\\
&&~~~~~\ldots,\tau_mE_m,0\ldots,0).
\end{eqnarray*}
Note that we have four kinds of block pairs $(A_i, D_i)$: $(\tau_iE_i,\tau_iE_iJ(\kappa_i,n_i))$, $(\tau_iE_i,\tau_iF_i)$, $(\tau_iF_i,\tau_iE_i)$ and $(0,0)$.
In fact, the second kind of block pairs is a special case of the first kind with $\kappa_i=0$ due to $E_iJ(0,n_i)=F_i$. We call the first two  kinds of block pairs type A block pairs, the third kind of block pairs type B block pairs and the last one type C block pairs.

\subsection{ SOCP reformulation from canonical form}
Without loss of generality, we make the following assumptions.
\begin{asmp}\label{A2.1}
i) There is at least one feasible solution in problem ${\rm (P)}$;
ii) The following three conditions do not hold true at the same time: $A\succeq0$, $b\in {\rm Range}(A)$ and $c=\frac{1}{2}b^TA^{+}b$.
\end{asmp}
Note that problem ${\rm (P)}$ is infeasible if and only if   $A\succeq0$, $b\in {\rm Range}(A)$ and $c=\frac{1}{2}b^TA^{+}b+k$ for some $k>0$, which leads to $h(x)=\frac{1}{2}(x+A^+b)^TA(x+A^+b)+k>0$.
If $A\succeq0$, $b\in {\rm Range}(A)$ and $c=\frac{1}{2}b^TA^{+}b$, then $h(x)=\frac{1}{2}x^TAx+b^Tx+c=\frac{1}{2}(x+A^+b)^TA(x+A^+b)\ge0$.
Thus, the inequality constraint in problem ${\rm(P)}$ becomes an equality constraint which means all the feasible solutions are in the boundary.
Actually, Assumption 2.1 is equivalent to Slation condition, i.e., there exists an $\overline{x}$ Such that $h(\overline{x})<0$.

Moreover, when the three conditions in ii) hold together, problem (P) reduces to an unconstrained quadratic problem: Decompose $A$ as $A=L^TL$, where $L\in \Re^{r\times n}$ with $r$
being the rank of $A$. Then, the constraint becomes $(x+A^+b)^TL^TL(x+A^+b)=0$ and thus $L(x+A^+b)=0$.
Rewrite $x=-A^+b+Vy$, where $V\in\Re^{n\times(n-r)}$ is a matrix basis of the null space  of $L$ and $y\in\Re^{n-r}$.
The problem then reduces to an unconstrained quadratic optimization problem.

Let us recall the S-lemma \cite{polik2007survey}, which states the equivalence of the following two statements under Slation condition:\\
$~~~~{\rm(S_1)}$ $(\forall x\in \Re^n)$ $h(x)\leq0$ $\Rightarrow$ $f(x)\geq0$.\\
$~~~~{\rm(S_2)}$ $\exists$ $\mu\geq 0$ such that $f(x)+\mu h(x)\geq0, ~\forall x\in\Re^n$.\\
In general, $f(x)$ can be represented as $f(x)=\frac{1}{2}x^TDx+e^Tx+v$ with an additional constant $v$ in $\rm (S_1)$ and $\rm (S_2)$, and we use this representation of $f(x)$ in the following of this section
when discussing about the S-lemma and its variants.
The connection between problem (P) and the S-Lemma is illustrated in \cite{xia2014s} by
\begin{eqnarray}
&v(P)&=\inf_{x\in \Re^n}\{~ f(x)\mid h(x)\leq 0\}\notag\\
&&=\sup_{\eta\in \Re}\{~\eta:\{x\in \Re^n\mid f(x)< \eta,h(x)\leq0\}=\emptyset\}\notag\\
&&=\sup_{\eta\in \Re}\{~\eta\mid \exists \mu\geq0 \text{ such that } f(x)-\eta+\mu h(x)\geq0,~\forall x\in \Re^n\}\notag\\
&&=\sup_{\eta\in \Re,\mu\geq0}\{~\eta\mid\left(
\begin{array}{ccc}
D+\mu A&  e+\mu b \\
e^T+\mu b^T &2(v+\mu c-\eta)
\end{array}
\right)\succeq0\}.\label{S2}
\end{eqnarray}

By invoking the S-lemma, we have the following theorems.
\begin{thm}\label{T2.2}
Consider the case where a type A block pair $(\tau_iE_i,\tau_iE_iJ(\kappa_i,n_i))$ exists in problem {\rm${\rm(P)}$}. If the size of the associated Jordan block $J(\kappa_i,n_i)$ is greater than $2$ and the associated eigenvalue of the Jordan block is real, then the objective value of ${\rm(P)}$ is unbounded from below, i.e., $v{\rm (P)}=-\infty$.
\end{thm}
\proof
If the size of the associated  Jordan block $J(\kappa_i,n_i)$  is greater than $2$, then $\tau_i(E_iJ(\kappa_i,n_i)+\mu E_i)$ takes the following form
 $$\tau_i(E_iJ(\kappa_i,n_i)+\mu E_i)=\tau_i\left(\begin{array}{ccccc}
& & &\kappa_i+\mu\\
& &\kappa_i+\mu&1 \\
&\cdots&\cdots&\\
\kappa_i+\mu&1 &&&
\end{array}\right).$$
Since the $(n_i-1)\times (n_i-1)$ principal minor $$\tau_i\left(
\begin{array}{ccccc}
&&\kappa_i+\mu&1\\
&\cdots&\cdots&\\
\kappa_i+\mu&1&\\
1&&&\end{array}
\right)$$
is non-positive semi-definite when its size $n_i-1$ is greater than or equal to 2,  $D_i+\mu A_i=\tau_i(E_iJ(\kappa_i,n_i)+\mu E_i)$ cannot be positive semi-definite.
Thus,  there is no $\mu$ such that $D+\mu A={\rm diag}(D_1+\mu A_1,\ldots,D_m+\mu A_m,0,\ldots,0)\succeq 0$. So the problem in (\ref{S2}) is infeasible and  by the S-lemma we have $v{\rm(P)}=-\infty$.\endproof

Using similar proofs, we have the following theorems.
\begin{thm}\label{T2.3}
Consider the case where a type A block pair $(\tau_iE_i,\tau_iE_iJ(\kappa_i,n_i))$ exists in problem {\rm${\rm(P)}$}. If the eigenvalues of the associated  Jordan block $J(\kappa_i,n_i)$ form a complex pair, then the objective value of  problem $\rm {(P)}$ is unbounded from below, i.e., $v{\rm (P)}=-\infty$.
\end{thm}
\proof
If the eigenvalues of the associated  Jordan block $J(\kappa_i,n_i)$ form a complex pair, then there does not exist a $\mu\in\Re$ such
that
$$\tau_i(E_iJ(\kappa_i,n_i)+\mu E_i)=\tau_i\left(\begin{array}{ccccc}
 & & &b_i &a_{i}+\mu\\
 & & &a_{i}+\mu&-b_{i} \\
&&\ldots&&\\
b_i &a_{i}+\mu&&&\\
a_{i}+\mu&-b_{i} &&&
\end{array}\right)\succeq0,$$ because the $4\times4$ principal minor (if $n_i=4k$ for some positive integer $k$),$$\tau_i\left(\begin{array}{ccccc}
 & & b_i &a_{i}+\mu\\
 & & a_{i}+\mu&-b_{i} \\
b_i &a_{i}+\mu&&\\
a_{i}+\mu&-b_{i} &&
\end{array}\right)$$
or the $2\times2$ principal minor (if  $n_i=4k+2$ for some positive integer $k$), $$\tau_i\left(
\begin{array}{cc}
b_i&a_{i}+\mu\\
a_{i}+\mu&-b_i
\end{array}
\right)$$
is non-positive semi-definite.
So the problem in (\ref{S2}) is infeasible and  by the S-lemma we get $v{\rm(P)}=-\infty$.
\endproof

\begin{thm}\label{T2.4}
Consider the case where a type B block pair $(\tau_iF_i,\tau_iE_i)$ exists in problem ${\rm(P)}$. If $\dim F_i\geq2$, then problem $\rm {(P)}$ is unbounded from below, i.e., $v{\rm (P)}=-\infty$.
\end{thm}
\proof
If the size of the associated  Jordan block $J(\kappa_i,n_i)$ is larger than or equal to 2, then there does not exist a $\mu\in\Re$ such
that
$$\tau_i(E_i+\mu F_i)=\tau_i\left(\begin{array}{ccccc}
 & & & &1\\
 & & &1&\mu \\
&&\ldots&&\\
 &1&\mu&&\\
1&\mu&&&
\end{array}\right)\succeq0.$$
So the problem in (\ref{S2}) is infeasible and  by the S-lemma we have $v{\rm(P)}=-\infty$.
\endproof

\begin{rem}
Assumption \ref{L2.1} is necessary in Theorem \ref{T2.4}. Otherwise the S-lemma does not hold and we have the following  counter example: $\min f(x)=x_1x_2 $ subject to $h(x)=x_2^2\le0$. The problem has a size $2$ block pair $(F_{2\times2}, E_{2\times2})$ but a finite optimal value of $\min\limits_{h(x)\le0} f(x)=0$.
\end{rem}

So if problem ${\rm (P)}$ has a finite optimal solution, then any type B block pairs are of size $1$ and any type A block pairs are of a size less than or equal to $2$ and the eigenvalues in the associated Jordan blocks are real. Now let us consider a type A block pair with size $2$, and, without loss of generality, let it be the first block $(A_1,D_1)=(\tau_1E_1,\tau_1E_1J_{1}(\lambda,2))$ with
\begin{eqnarray*}
E_1=
\left(
\begin{array}{cc}
0 &  1 \\
1& 0
\end{array}
\right),~
E_1J_1(\lambda,2)=
\left(
\begin{array}{ccc}
0  &  \lambda \\
\lambda & 1
\end{array}
\right).
\end{eqnarray*}
Denote $D= {\rm diag}(D_1,D_{J})$, $A= {\rm diag}(A_1,A_J)$, $I=\{1,2\},J=\{3,\ldots,n\}$, $e_I=(e_1,e_2)^T$,
 $e_J=(e_{3},\ldots,e_n)^T$, $b_I=(b_1,b_2)^T$, $b_J=(b_{3},\ldots,b_n)^T$, $z=(x_1,x_2)^T$, and $y=(x_{3},\ldots,x_n)^T$. We can then represent problem ${\rm (P)}$ as follows,
\begin{eqnarray*}
{\rm (RP)}~&\text{min}&\frac{1}{2}y^TD_Jy+e_J^Ty+\frac{1}{2}z^TD_1z+e_I^Tz \\
          &{\rm s.t.}&\frac{1}{2}y^TA_Jy+b_J^Ty+\frac{1}{2}z^TA_1z+b_I^Tz+c\leq 0.\\
\end{eqnarray*}
The term in the constraint associated with $(A_1,D_1)$ is
\begin{equation}
\frac{1}{2}z^TA_1z+b_I^Tz=\frac{1}{2}\tau_1z^TE_1z+b_I^Tz=\tau_1z_1z_2+b_1z_1+b_2z_2,\label{cons}
\end{equation}
and the term in the objective function associated with $(A_1,D_1)$ is
\begin{equation}
\frac{1}{2}z^TD_1z+e_I^Tz=\frac{1}{2}\tau_1z^TE_1J_1(\lambda,2)z+e_I^Tz=\tau_1\lambda z_1z_2+\frac{1}{2}\tau_1z_2^2+e_1z_1+e_2z_2.\label{obj}
\end{equation}
Without loss of generality, we further assume $b_1=b_2=0$. Since otherwise when letting $z_1'=z_1+\tau_1b_2$ and $z_2'=z_2+\tau_1b_1$,
the constraint function will become $\frac{1}{2}y^TA_Jy+b_J^Ty+\frac{1}{2}{z'}^TA_1{z'}+c'$, where $c'=c-\tau_{1} b_1b_2$,
and the objective function will become $\frac{1}{2}y^TD_Jy+e_J^Ty+\frac{1}{2}{z'}^TD_1{z'}+{e'}_I^T{z'}+d_0$, where
${e'}_1=e_1-\lambda b_1$, ${e'}_2=e_2-b_{1}-\lambda b_2$, and $d_0=-e_1\tau_1b_2-e_2\tau_1b_1+\tau_1\lambda b_1b_2+\frac{1}{2}\tau_1 b_1^2$.
Note that $\tau_1=\pm 1$ according to Theorem \ref{T2.1}.

From now on, we assume that the coefficients in $b$ corresponding to any $2\times 2$ type A Jordan block pair are $0$.

\begin{thm}\label{T2.5}
Consider the case where there exists a type A block pair $(\tau_1E_1,\\\tau_1E_1J_{1}(\lambda,2))$ in problem {\rm (P)} and the eigenvalue of the associated  Jordan block $J_{1}(\lambda,2)$ is  real. Assume there is a feasible solution $\bar{x}=(\bar{z}^T,\bar{y}^T)^T$ and let
$\pi=\tau_1\bar{z}_1\bar{z}_2$.
Let $\rho=\inf \{~(\ref{obj})  \mid(\ref{cons})=\tau_1z_1z_2\leq\pi\}$. We have the following three cases:
\begin{enumerate}
\item When $\tau_1=1$. If $(\lambda\leq0,~e_1=0,~e_2\neq0)$ or $(\lambda=0,~e_1=0,~e_2=0,~\pi\geq0)$ or $(\lambda<0,~e_1=0,~ e_2=0, ~\pi=0)$, then $\rho=\lambda\pi-\frac{1}{2}e_2^2$ and the infimum is attainable;
\item  When $\tau_1=1$. If $(\lambda=0,~e_1=0,~ e_2=0,~ \pi<0)$ or $(\lambda<0,~ e_1=0,~ e_2=0,~ \pi\neq0)$ , then $\rho=\lambda\pi-\frac{1}{2}e_2^2$ and the infimum is unattainable;
\item Otherwise, $\rho=-\infty$  and thus  problem ${\rm(P)}$ is unbounded from below.
\end{enumerate}
\end{thm}
\proof
We consider the problem in the following cases:
\begin{itemize}
\item When $\tau_1=1$, $(\ref{cons})$ becomes $z_1z_2\leq\pi$ and (\ref{obj}) becomes $\lambda z_1z_2+\frac{1}{2}z_2^2+e_1z_1+e_2z_2$. We then have the following cases corresponding to the values of $\lambda$, $e_1$, $e_2$ and $\pi$.
\begin{itemize}
\item
When $\lambda>0$, set $z_1=-\frac{M}{\lambda}-M$ and $z_2=M$, where $M\in\Re$ is chosen such that $-(1+\frac{1}{\lambda})M^2\leq \pi$. Then $z_1z_2\leq\pi$ and $\rho=-(\lambda+\frac{1}{2}) M^2-(\frac{e_1}{\lambda}+ e_1 -e_2) M\rightarrow-\infty$ when $M\rightarrow\infty$. (case 3)

\item When $\lambda=0$, we have the following subcases:
\begin{itemize}
\item When $ e_1\neq0$, set $ z_1=-\text{sign}(e_1)M$ and $z_2=-\text{sign}(e_1)\frac{\pi}{M}$, where $M\in\Re$ and $z_1z_2=\pi$.
Then $\rho=-|e_1|M+C_1+C_2\cdot \frac{1}{M}+C_3\cdot \frac{1}{M^2}\rightarrow-\infty$ when $M\rightarrow+\infty$, where $C_i, i=1, 2, 3$, are the reduced constants. (case 3)
\item When $e_1=0$, (\ref{obj}) becomes $\frac{1}{2}z_2^2+e_2z_2=\frac{1}{2}(z_2+e_2)^2-\frac{1}{2}e_2^2\geq-\frac{1}{2}e_2^2$ $\Rightarrow$ $\rho\geq-\frac{1}{2}e_2^2$.
\begin{itemize}
\item If $e_2\neq0$, set $z_1=-\frac{\pi}{e_2}$ and $z_2=-e_2$, then $z_1z_2=\pi$ and $(\ref{obj})=-\frac{1}{2}e_2^2$. Thus $\rho=-\frac{1}{2}e_2^2$ and the infimum is attainable. (case 1)
\item If $e_2=0$ and $\pi\geq0$, we can set $z_1$ at any real value and $z_2=0$ such that $z_1z_2\le\pi$ and thus $\rho=-\frac{1}{2}e_2^2$ and the infimum is attainable. (case 1)
\item If otherwise $e_2=0$ and $\pi<0$, we cannot set $z_2=-e_2=0$, which contradicts the constraint $z_1z_2\leq\pi<0$. So the infimum is unattainable.
 But we can set $z_1=M\pi$ and $z_2=\frac{1}{M}$ ($M\in\Re$) such that$z_1z_2=\pi$ and $(\ref{obj})=\frac{1}{2M^2}\rightarrow 0$ when $M\rightarrow\infty$. Thus $\rho=\lim\limits_{M\rightarrow \infty}\frac{1}{2M^2}-\frac{1}{2}e_2^2=0$ but the infimum is unattainable. (case 2)
\end{itemize}
\end{itemize}
\item When $\lambda<0$,  we have the following subcases:
\begin{itemize}
\item When $e_1\neq0$, set $z_1=-\text{sign}(e_1)M$ and $z_2=-\frac{\text{sign}(e_1)\pi}{M}$ ($M\in\Re$) such that $z_1z_2=\pi$ and then $(\ref{obj})=-|e_1|M+C_1+C_2\cdot \frac{1}{M}+C_3\cdot \frac{1}{M^2}\rightarrow-\infty$ when $M\rightarrow\infty$, where $C_i,~ i=1,2,3$ are the reduced constants.
Thus $\rho=-\infty$. (case 3)
\item When $e_1=0$, $(\ref{obj})=\lambda z_1z_2+\frac{1}{2}z_2^2+e_2z_2\geq\lambda\pi+\frac{1}{2}(z_2+e_2)^2-\frac{1}{2}e_2^2\geq\lambda\pi-\frac{1}{2}e_2^2$. Next we show that $\rho=\lambda\pi-\frac{1}{2}e_2^2$. We first note that, to achieve $\lambda z_1 z_2$ = $\lambda \pi$ in the above inequality, we need to set $z_1 z_2$ = $\pi$.
\begin{itemize}
\item If $e_2\neq0$, set $z_1=-\frac{\pi}{e_2}$ and $z_2=-e_2$, such that $z_1z_2=\pi$ and then $(\ref{obj})=\lambda \pi-\frac{1}{2}e_2^2$. (case 1)
\item If $e_2=0$ and $\pi\neq0$, we cannot set $z_2=-e_2=0$, which contradicts the constraint $z_1z_2 = \pi \neq 0$. So the infimum is unattainable. But we can set $z_1=M\pi$ and $z_2=\frac{1}{M}$ ($M\in\Re$) such that $z_1z_2=\pi$ and $(\ref{obj})=\lambda \pi+\frac{1}{2M^2}\rightarrow \lambda \pi$ when $M\rightarrow\infty$. So $\rho=\lambda \pi-\frac{1}{2}e_2^2=\lambda \pi$  and the infimum is unattainable. (case 2)
\item If $e_2=0$ and $\pi=0$, we can set $z_1$ at any real value and $z_2=0$ and thus attain the infimum $\rho=\lambda \pi-\frac{1}{2}e_2^2=\lambda \pi$. (case 1)
\end{itemize}
\end{itemize}
\end{itemize}
\item When $\tau_1=-1$, set $z_1=-\frac{\pi}{M}$ and $z_2=M$ ($M\in\Re$) such that $\tau_1z_1z_2=-z_1z_2=\pi$ and $(\ref{obj})=-\frac{1}{2}M^2+C_1\cdot M+C_2\cdot \frac{1}{M}+C_3\rightarrow-\infty$ when $M\rightarrow\infty$, where $C_i,~ i=1, 2, 3$, are the reduced constants.
Thus $\rho=-\infty$. (case 3)
\end{itemize}
Since $\inf \{~(\ref{obj})  \mid(\ref{cons})=\tau_1z_1z_2\leq\pi\}$ is a subproblem of {\rm(RP)}, if there is a feasible solution $(\bar{z}^T,\bar{y}^T)^T$ for {\rm(RP)} with $\bar{z}_1\bar{z}_2=\pi$, and $\rho=-\infty$, then $v{\rm (P)}=v{\rm (RP)}=-\infty$.
\endproof

\begin{rem}\label{rem2.3}
\begin{enumerate}
\item Case 2 in the above theorem is the only case where problem $\rm{\rm{(P)}}$ is bounded from below but its infimum is unattainable.
\item If two matrices $A$ and $D$ are simultaneously diagonalizable via congruence, which covers most conditions discussed in the existing literature (\cite{ben2014hidden,feng2012duality,more1993generalizations}), then the optimal value is attainable when problem ${\rm (P)}$ is bounded from below.
\item
All the bounded cases require $e_1=0$ in the linear terms in the objective function associated with $2\times 2$ Jordan blocks.
\end{enumerate}
\end{rem}

The following theorem is a direct result of Theorems \ref{T2.2}--\ref{T2.5}.
\begin{thm}\label{C2.1}If problem $\rm(P)$ has an optimal value bounded from below, then:
\begin{enumerate}
\item $\dim{E_i}\leq2,~ i=1,\ldots,p$, $\dim{E_i}=1,~i=p+1,\ldots,m,$ and there is no complex eigenvalue pair in $J(\kappa_i,n_i)$;
\item If for some index $i$, $\dim{E_i}=2$, then the $i$th block satisfies the conditions in either case 1 or case 2 in Theorem \ref{T2.5}.
\end{enumerate}
\end{thm}

Note that the conditions in items $1$ and $2$ of Theorem \ref{C2.1}  are necessary for problem (P) to be bounded from below and we assume these conditions
hold in the following discussion of this section.  Rearrange the block pairs with single elements to the upper left part of the diagonal
in the canonical form and express $A$ and $D$ in the following forms, \begin{equation} A={\rm diag}(\alpha_1,\ldots,\alpha_l,E_1,\ldots,E_{\frac{n-l}{2}}),\label{FinA}\end{equation}
\begin{equation}
D={\rm diag}(\delta_1,\ldots,\delta_l,E_1J(\zeta_1,2),\ldots,E_{\frac{n-l}{2}}J(\zeta_{\frac{n-l}{2}},2)),\label{FinB}
\end{equation} where $l$ is the number of block pairs with size $1$ in the canonical form of (\ref{CF1}) and (\ref{CF2}).
In the following, we assume $A$ and $D$ are in the form of $(\ref{FinA})$ and $(\ref{FinB})$  and further assume $b_i=0,~i = l+1, \ldots, n$, as we discussed after (\ref{obj}).

Moreover, from item 3 in Remark \ref{rem2.3}, we have  $e_{l+2j-1}=0,~j=1,\ldots,\frac{n-l}{2}.$
Then  problem   ${\rm (P)}$ can be reduced to the following form:
\begin{eqnarray*}
{\rm (P_1)}~&\text{min}& f(x)=\sum_{i=1}^l(\delta_i x_i^2+e_i x_i)+\sum_{j=1,\ldots,\frac{n-l}{2}}(\zeta_jx_{l+2j-1}x_{l+2j}+\frac{1}{2}x_{l+2j}^2+e_{l+2j}x_{l+2j})\\
&{\rm s.t.}&h(x)=\sum_{i=1}^l(\alpha_i x_i^2+b_i x_i)+\sum_{j=1,\ldots,\frac{n-l}{2}}(x_{l+2j-1}x_{l+2j})+c\leq0.
\end{eqnarray*}

\begin{thm}\label{mainthm}
Assume that items 1 and 2 in Theorem \ref{C2.1} are satisfied, then $v{\rm(P)}=v{\rm(P_1)}=v{\rm(P_2)}$, where ${\rm(P_2)}$ is the following SOCP problem:
\begin{eqnarray*}
{\rm (P_2)}~&~\min &\sum_{i=1}^l(\delta_i y_i+e_i x_i)+\sum_{j=1}^\frac{n-l}{2}\zeta_j z_j+c_0 \\
&~{\rm{s.t.}}&\sum_{i=1}^l(\alpha_i y_i+b_i x_i)+\sum_{j=1}^\frac{n-l}{2}z_j+c\leq0,\\
&~&\frac{1}{2}x_i^2-y_i\leq0,~\forall i=1,2,\ldots,l,\\
&~&x,y\in \Re^l, ~z\in \Re^{\frac{n-l}{2}},
\end{eqnarray*}
where $c_0=-\sum_{j=1,\ldots,\frac{n-l}{2}}\frac{1}{2}e_{l+2j}^2$.

More specifically, if ${\rm(P_2)}$  admits an optimal solution, then there exists an optimal solution  $(\bar{x},\bar{y},\bar{z})$ to ${\rm(P_2)}$
 with $\frac{1}{2}\bar{x}_i^2=\bar{y}_i,~i=1,2,\ldots,l$.
Moreover, we can find  an optimal solution (or an $\epsilon$ optimal solution)   $\tilde{x}$ to ${\rm(P_1)}$ with
\begin{eqnarray}
\left.\begin{array}{ll}
\tilde{x}_i=\bar{x}_i,i=1,\ldots,l,\\
\tilde{x}_{l+2j}=\left\{\begin{array}{ll}
1/M&\text{if }\left\{\begin{array}{ll} \zeta_j=0,e_{l+2j}=0,\bar{z}_j<0,\\\text{or } \zeta_j<0,e_{l+2j}=0,\bar{z}_j \neq0,\end{array}\right.\\
-e_{l+2j}&\text{otherwise},
\end{array}\right. ~~~j=1,\ldots,\frac{n-l}{2},\\
\tilde{ x}_{l+2j-1}=\frac{\bar{z}_j}{\tilde{x}_{l+2j}},~j=1,\ldots,\frac{n-l}{2}.
\end{array}\right. \label{un1}
\end{eqnarray}

Particularly, if  ${\rm(P_1)}$ is bounded from below, then the optimal value of ${\rm(P_1)}$ is unattainable  if and only if $\zeta_j=0,~e_{l+2j}=0,~\bar{z}_j<0 \text{ or } \zeta_j<0,~e_{l+2j}=0,~\bar{z}_j\neq0$.
In this case, for any $\epsilon>0$, there exists an $\epsilon$ optimal solution  $\tilde{x}$ such that $f(\tilde{x})-v{\rm(P_1)}<\epsilon$ with a sufficiently large $M>0$.
 \end{thm}

\proof
Because of Theorem \ref{C2.1},  ${\rm(P_1)}$ is equivalent to ${\rm(P)}$. And the main differences between ${\rm(P_1)}$ and ${\rm(P_2)}$ are the terms associated
to the $2\times2$ Jordan blocks. Let us consider how to simplify the terms associated with the $2\times2$ Jordan blocks.
According to Assumption \ref{A2.1}, ${\rm(P_1)}$ is feasible. For any feasible solution $\hat{x}$ of ${\rm(P_1)}$, we let $\pi_{j}= \hat{x}_{l+2j-1}\hat{x}_{l+2j}$. Now let us concentrate on problem $\inf\{\zeta_jx_{l+2j-1}x_{l+2j}+\frac{1}{2}x_{l+2j}^2+e_{l+2j}x_{l+2j}\mid x_{l+2j-1}x_{l+2j}=\pi_j\}$ = $\inf\{\zeta_j\pi_j+\frac{1}{2}x_{l+2j}^2+e_{l+2j}x_{l+2j}\mid x_{l+2j-1}x_{l+2j}=\pi_j\}$ = $\inf\{\zeta_j\pi_j+\frac{1}{2}(x_{l+2j}^2+e_{l+2j})^2 - \frac{1}{2}e^2_{l+2j}\mid x_{l+2j-1}x_{l+2j}=\pi_j\}$.
Thus setting $x_{l+2j}=-e_{l+2j}$ (if $e_{l+2j}=0$, set $x_{l+2j}=\frac{1}{M}$ as in the proof of Theorem \ref{T2.5})
and $x_{l+2j-1}=\frac{\pi_{j}}{x_{l+2j}}$, the objective function $(\zeta_{j}x_{l+2j-1}x_{l+2j}+\frac{1}{2}x_{l+2j}^2+e_{l+2j}x_{l+2j})$ has an infimum $\zeta_j\pi_{j}-\frac{1}{2}e_{l+2j}^2$ under the constraint $x_{l+2j-1}x_{l+2j}=\pi_j$, which is linear with the cross term $x_{l+2j-1}x_{l+2j}=\pi_{j}$.

Using such a separability, we denote $$z_j=x_{l+2j-1}x_{l+2j} {\rm~ and~} c_0=-\sum_{j=1,\ldots,\frac{n-l}{2}}\frac{1}{2}e_{l+2j}^2,$$ and have the following problem
which has the same objective value with ${\rm(P_1)}$:
\begin{eqnarray*}
{\rm (P_3)}&\min &\sum_{i=1}^l(\delta_i x_i^2+e_i x_i)+\sum_{j=1}^\frac{n-l}{2}\zeta_j z_j+c_0 \\
&{\rm{s.t.}}&\sum_{i=1}^l(\alpha_i x_i^2+b_i^T x_i)+\sum_{j=1}^\frac{n-l}{2} z_j+c\leq0.
\end{eqnarray*}
Moreover, if there is an optimal solution $(\bar{x}, \bar{z})$ of ${\rm(P_3)}$,
we can also find  an optimal solution (or an $\epsilon$ optimal solution)  $\tilde{x}$ of ${\rm(P_1)}$ in the  form of $(\ref{un1})$.
In this case, the optimal value of ${\rm(P_1)}$ is unattainable if and only if $\zeta_j=0,~e_{l+2j}=0,\bar{z}_j<0 \text{ or } \zeta_j<0,~e_{l+2j}=0,~\bar{z}_j\neq0$  from Theorem \ref{T2.5}. Furthermore, for any $\epsilon>0$, if we set $M\geq\sqrt \frac{1}{2\epsilon}$,  then  $f(\tilde{x})-v({\rm P_1})=\frac{1}{2M^2}\leq\epsilon$.

Introducing $y_i=\frac{1}{2}x_i^2,~ i=1,2,\ldots,l$, ${\rm (P_3)}$ is then equivalent to the following ${\rm (P_4)}$:
\begin{eqnarray*}
{\rm (P_4)}~&\min&\sum_{i=1}^l(\delta_i y_i+e_i x_i)+\sum_{j=1}^\frac{n-l}{2}\zeta_j z_j+c_0 \\
&{\rm{s.t.}}&\sum_{i=1}^l(\alpha_i y_i+b_i x_i)+\sum_{j=1}^\frac{n-l}{2}z_j+c\leq0\\
&~&\frac{1}{2}x_i^2-y_i=0,~\forall i=1,2,\ldots,l,\\
&~&x,y\in \Re^l, ~z\in \Re^{\frac{n-l}{2}}.
\end{eqnarray*}

We next prove the equivalence of $\rm{(P_4)}$ and $\rm{(P_2)}$ in two parts:
\begin{enumerate}
\item If ${\rm(P_2)}$ is unbounded from below, then ${\rm(P_4)}$ is unbounded from below.
\item If $({\rm P_2})$ has an optimal solution $(x^*,y^*, z^*)$, then we can always find a solution $(\bar{x},\bar{y},\bar{z})$ with
$\bar{y}_i=\frac{1}{2}\bar{x}_i^2,~i=1,\ldots,l$ and $\bar{z}=z^*$, which is  optimal not only to ${\rm (P_2)}$ but also to ${\rm (P_4)}$.
\end{enumerate}

The first part is proved in the following Lemma \ref{L2.3}.
Now let us prove part 2.
Note that if ${\rm(P_2)}$ is bounded from below, then there must exist an optimal solution $(x^*,y^*, z^*)$ since Slation condition is satisfied.   Denote
\begin{eqnarray*}
J:=\{i:\frac{1}{2}(x_i^*)^2<y_i^*,~i=1,\ldots,l\}.
\end{eqnarray*}
If $J=\emptyset$, then $(x^*,y^*,z^*)$ is also an optimal solution of $\rm{(P_4)}$.
If $J\neq\emptyset$, by Theorem 3 in \cite{ben2014hidden}, we can transform the optimal solution $(x^*,y^*,z^*)$  of $\rm{(P_2)}$ to an optimal solution $(\bar{x},\bar{y},z^{*})$  of $\rm{(P_2)}$ with $\bar{y}_i=\frac{1}{2}\bar{x}_i^2$, $i=1,\ldots,l$, and $(\bar{x},\bar{y},z^{*})$ is also a feasible solution of ${\rm (P_4)}$, since  $\bar{y}_i=\frac{1}{2}\bar{x}_i^2$, $i=1,\ldots,l$.
So $v({\rm{P_2}})\geq v({\rm{P_4}})$. But ${\rm (P_2)}$ is a relaxation of ${\rm (P_4)}$, so $v({\rm P_2})\leq v({\rm P_4})$. Thus $v({\rm P_2})=v({\rm P_4})$ and $(\bar{x},\bar{y},z^{*})$ is  optimal to ${\rm (P_4)}$.
\endproof

\begin{lem}\label{L2.3}
If ${\rm(P_2)}$ is unbounded from below, then ${\rm(P_4)}$ is unbounded from below.
\end{lem}
\proof
We only need to prove that ${\rm(P_3)}$ (since ${\rm(P_4)}$ is equivalent to ${\rm(P_3)}$) is bounded from below implies that ${\rm(P_2)}$ is bounded from below.

In this proof, we only consider the cases with no $z$ term in ${\rm(P_2)}$ and ${\rm(P_3)}$, since $z$ only appears in linear terms in both the objective and constraint functions., which can be regarded as a special case of the $x$ variable (i.e., the coefficients before $z_j^2$ are $0$, $j=1,\ldots,\frac{n-l}{2}$).

Denote the Lagrangian function of ${\rm(P_3)}$ as $L(x,\nu)=f(x)+\nu h(x)$, and the dual function as $\theta(\nu)=\min\limits_x L(x,\nu)$, where $\nu\geq0$.
If ${\rm(P_3)}$ is bounded from below, from the S-lemma (as Slation condition holds here), we know there is no duality gap between the primal problem ${\rm(P_3)}$
and its Lagrangian dual problem of $\max\limits_{\nu\geq0} \theta(\nu)$, i.e., there exists $(\bar{x},\bar{\nu})$ such that $\min \limits_{h(x)\leq0} f(x)=f(\bar{x})=\theta(\bar{\nu})=\max\limits_{\nu\geq0}\theta(\nu)$.
So $(\bar{x},\bar{\nu})$ is a saddle point of the Lagrangian function $L(x,\nu)$.
Then $f(\bar{x})=  \min\limits_{x\in {\Re^n}}  L(x,\bar{\nu})$,  $h(\bar{x})\leq0$, $\bar{\nu}\geq0$, $\bar{\nu} h(\bar{x})=0$.
From  $\min\limits_{x\in {\Re^n}} L(x,\bar{\nu})=\min\limits_{x\in {\Re^n}}\sum_{i=1}^l(\frac{1}{2}(\delta_i+\bar{\nu} \alpha_i)x_{i}^2+(e_i+\bar{\nu} b_i)x_{i})+c_{0}+\bar{\nu} c=f(\bar{x})$,
we get $\delta_i+\bar{\nu} \alpha_i\geq0$ and $(\delta_i+\bar{\nu} \alpha_i)\bar{x}_{i}+(e_i+\bar{\nu} b_i)=0$ and  if, in addition, $\delta_i+\bar{\nu} \alpha_i=0$, we have $e_i+\bar{\nu} b_i=0$, $i=1,\ldots,l$.
So $(\bar{x},\bar{\nu})$ satisfies the KKT conditions of ${\rm(P_3)}$ , i.e.,  $(\delta_i+\bar{\nu} \alpha_i)\bar{x}_i+e_i+\bar{\nu} b_i=0$, $i=1,\ldots,l$, $\bar{\nu}\geq0$,
 $\bar{\nu}h(\bar{x})=0$.

Next we can construct a KKT point of ${\rm(P_2)}$ from the saddle point $(\bar{x},\bar{\nu})$. Denote $\bar{y}_i=\frac{1}{2}\bar{x}_i^2$ and $\bar{\mu}_i= \delta_i+\bar{\nu} \alpha_i\geq0$, $i=1,\ldots,l$. Then  $(\bar{x},\bar{y},\bar{\mu},\bar{\nu})$ satisfies the KKT condition of ${\rm(P_2)}$:
$\delta_i+\bar{\nu} \alpha_i-\bar{\mu}_i=0,e_i+\bar{\nu} b_i+\bar{\mu}_i\bar{x}_i=0$, $\bar{\mu}_{i}(\frac{1}{2}\bar{x}_i^2-\bar{y}_i)=0$, $i=1,\ldots,l$, $\bar{\nu}(\sum_{i=1}^l(\alpha_i\bar{y}_i+b_i\bar{x}_i)+c)=0$. Thus  $(\bar{x},\bar{y})$ is a global optimal solution of ${\rm(P_2)}$ because of the convexity of ${\rm(P_2)}$. So we conclude that  ${\rm(P_2)}$ is bounded from below.
\endproof

\begin{exam} Consider the following problem:
\begin{eqnarray*}
~&~\min& -x_1x_2+0.5x_2^2-x_3^2+x_4^2+2x_2-x_4\\
&~{\rm s.t.}& x_1x_2+x_3^2+0.75x_4^2\leq1.25,
\end{eqnarray*}
where the related matrices can be expressed as
$$A=\left(
\begin{array}{cccc}
0 &  1 &0&0\\
1& 0&0&0\\
0&0&2&0\\
0&0&0&1.5
\end{array}
\right),~~D=\left(
\begin{array}{cccc}
0&-1&0&0\\
-1&1&0&0\\
0&0&-2&0\\
0&0&0&2
\end{array}
\right),$$
$e=(0,2,0,-1)^T$, $b=0$ and $c=-1.25$. Note that $A$ and $D$ are not simultaneously diagonalizable but in the canonical form (\ref{CF1}) and (\ref{CF2}).
According to Theorem \ref{mainthm}, we get the following equivalent SOCP reformulation,
\begin{eqnarray*}
~&~\min&z_1-2y_3+2y_4-x_4-2\\
&~{\rm s.t.}&z_1+2y_3+1.5y_4\leq1.25\\
&~&\frac{1}{2}x_3^2-y_3\leq0\\
&~&\frac{1}{2}x_4^2-y_4\leq0.
\end{eqnarray*}
Solving the above SOCP problem yields the optimal solution $z_1^*=-12.9773
$, $x_3^*=0$, $x_4^*=0.2857$, $y_3^*=7.0830$ and $y_4^*=0.0408$. Note  $\frac{1}{2}(x_3^*)^2-y_3^*=-7.0830<0$. Using the transformation method in Theorem 3 in \cite{ben2014hidden}, we obtain a new solution $(\bar{x},\bar{y},\bar{z})$, with $\bar{x}_3=\pm\sqrt{2y_3^{*}}=\pm3.7638$, $\bar{z}_1=z_1^*$, $\bar{x}_3=x_3^*$, $\bar{y}_3=y_3^*$ and $\bar{y}_4=y_4^*$.
By applying Theorem \ref{mainthm},  we get $\tilde{x}_2=-2$, $\tilde{x}_1=\frac{\bar{z}_1}{\tilde{x}_2}=6.4886$, $\tilde{x}_3=\bar{x}_3$ and $\tilde{x}_4=\bar{x}_4$. So we obtain an optimal solution  $\tilde{x}=(6.4886,-2,\pm3.7638,0.2857)^T$ to the origin problem, with an optimal value $   -3.3929$.
\end{exam}

\subsection{Dual problem of  SOCP reformulation}
The following theorem shows that the dual problem of $\rm{(P_2)}$ is  a simple concave maximization problem with a single variable.
\begin{thm}\label{T2.7}
Under Assumption 2.1, the objective values of $\rm{(P_2)}$ and the following Lagrangian dual
 problem are equal:\begin{eqnarray*}
({\rm{D_1}})~~
\max\limits_{\nu\geq0}\{\rho(\nu)=c\nu+c_0+\sum_{i=1}^l h_i(\nu)+g(\nu)\},
\end{eqnarray*}
where
\begin{eqnarray}
h_i(\nu)=
\left\{\begin{array}{ll}
-\frac{(\nu b_i+e_i)^2}{2(\nu\alpha_i+\delta_i)}&\text{if }\nu\alpha_i+\delta_i>0,\\
0&\text{if }\nu\alpha_i+\delta_i=0\text{ and }e_i+\nu b_i=0,\\
-\infty&\text{otherwise},
\end{array}\right. \label{h_expression}
\end{eqnarray}
$i=1,\ldots,l,$ and
\begin{eqnarray}
g(\nu)=\left\{\begin{array}{ll}
0&\text{if }\ \zeta_j+\nu=0,~\forall j=1,\ldots,\frac{n-l}{2},\\
-\infty&\text{otherwise}.
\end{array}\right. \label{g_expression}
\end{eqnarray}
\end{thm}
\proof

The Lagrangian function of ${\rm{(P_2)}}$ is:
\begin{eqnarray*}
\begin{array}{lll}
&L(x,y,z;\mu,\nu)&=\delta^T y+e^T x+\zeta^T z+c_0+\nu(\alpha^Ty+b^Tx+\bar{1}^Tz+c)\\
&&~~+\sum_{i=1}^l\mu_i(\frac{1}{2}x_i^2-y_i)\\
&&=\sum_{i=1}^ly_i(\delta_i+\nu\alpha_i-\mu_i)+\sum_{i=1}^lx_i(e_i+\frac{1}{2}\mu_ix_i+\nu b_i)\\
&&~~+\sum_{j=1}^{\frac{n-l}{2}}(\zeta_j+\nu) z_j+cv+c_0,
\end{array}
\end{eqnarray*}
where $\mu\geq0$ and $\nu\geq0$.

The dual objective function is then given as
\begin{eqnarray*}
m(\mu,\nu)=\min\limits_{x,y,z}L(x,y,z;\mu,\nu).
\end{eqnarray*}
If $m(\mu,\nu)$ is bounded from below, then the coefficients of $y_i$ and $z_i$ should be equal to $0$.
As the objective function is separable, we can get an analytical form of $m(\mu,\nu)$:
\begin{eqnarray*}
\max\limits_{\mu\geq0,\nu\geq0}m(\mu,\nu)=\max\limits_{\nu\geq0}(c\nu+c_0+\sum_{i=1}^lh_i(\nu)+g(\nu)),
\end{eqnarray*}
where $h_i(\nu)$ and $g(\nu)$ are given in (\ref{h_expression}) and (\ref{g_expression}), respectively.
\endproof

\begin{rem}\label{R2.4}
If  variable $z$  exists  in $\rm{(P_2)}$, i.e., $\rm{(P_2)}$ has some $2\times 2$ block pairs,
then either all $\zeta_i,~i=1,~\ldots,~\frac{n-l}{2}$, in $\rm{(P_2),}$ are the same (there exists a $\nu\geq0$ such that
$\nu+\zeta_i=0,~i=1,\ldots,\frac{n-l}{2}$, and $v{\rm(D_1)}=\rho(\nu)$),
or, otherwise, there does not exist such a $\nu$ and $v{\rm(P_2)}=v\rm(D_1)=-\infty.$ The first case requires that all $\zeta_i$ are equal,
which is equivalent to that all $2\times2$ Jordan blocks have the same eigenvalue, and  we can directly calculate the optimal value by the dual problem $(\rm D_1)$ using $\nu=\zeta_i$. $\forall i$.
In fact, we characterize all the possible situations of the GTRS: If the GTRS is bounded from below, then there
only exist $1\times1$ or $2\times2$ Jordan blocks in the canonical form and in addition all $2\times2$ Jordan blocks have the same eigenvalue.
While the results in Ben-Tal and Hertog \cite{ben2014hidden} are restricted to simultaneously diagonalizable cases, our results identify all the cases where the GTRS is bounded from below,
though only one situation (i.e., there exist $2\times2$ Jordan blocks in the canonical form and all $2\times2$ Jordan blocks have the same eigenvalue) is proved to be possibly bounded from below.

 If the variable $z$ does not exist in ${\rm(P_2)}$, ${\rm(D_1)}$ can be formulated as an SOCP:
\begin{eqnarray*}
\rm(D_1')~
&\min&c\nu+c_0+\sum_{i=1}^nw_i \\
&{\rm s.t.}&(\nu\alpha_i+\delta_i+\frac{1}{2}w_i)^2\geq(\nu\alpha_i+\delta_i-\frac{1}{2}w_i)^2+(\nu b_i+e_i)^2,~\forall i\\
&&\nu\alpha_i+\delta_i\geq0,~\forall i\\
&&\nu\geq0.
\end{eqnarray*}
\end{rem}

\subsection{A simplified S-lemma}
Denote $h(x)=\frac{1}{2}x^TAx+b^Tx+c$ and $f(x)=\frac{1}{2}x^TDx+e^Tx+v$ in the following of this section when discussing about the S-lemma and its variants.  First recall the classical S-lemma in Section 2.2 and note that LMI are involved in the statement {$\rm (S_2)$} from $(\ref{S2})$.

By applying the canonical form and the SOCP reformulation ${\rm(P_2)}$ of  problem ${\rm(P)}$, we can obtain  a  simplified S-lemma, which shows the equivalence of the following statements ${\rm(S_1)}$  and $\rm{(\tilde{S}_2)}$  when Slation condition holds:\\
$~~~~{\rm(S_1)}$ $(\forall x\in \Re^n)$ $h(x)\leq0\Rightarrow f(x)\geq0$.\\
$~~~~{\rm(\tilde{S}_2)}$ $\exists x,y\in \Re^l,~ z\in \Re^{\frac{n-l}{2}}:\left\{\begin{array}{ll}
\delta^T y+e^T x+\zeta^T z+c_0+v\geq0\\
\alpha^Ty+b^Tx+\bar{1}^Tz+c\leq0\\
\frac{1}{2}x_i^2-y_i\leq0,~\forall~ i=1,\ldots,l.\end{array}\right.$\\
When solving the SOCP problem in ${\rm(\tilde{S}_2)}$, we avoid LMI, which  are hard to handle for
large-scale problems.
 Moreover, our simplified version of S-lemma includes the simplified version of S-lemma in Ben-Tal and Hertog \cite{ben2014hidden} as a special case.

\subsection{Algorithms for computing the canonical form}
According to the proof of Theorem 1 in Uhlig \cite{uhlig1976canonical}, we present a method in the appendices to get the canonical form for two arbitrary $n\times n$ symmetric
matrices under the following conditions:
\begin{enumerate}
\item  $A$ is nonsingular;
\item
The Jordan normal form of $A^{-1}D$ has only $1\times1$ and $2\times2$ Jordan blocks with real eigenvalues and all the $2\times 2$ Jordan blocks have the same eigenvalue.\end{enumerate}
We assume that the first condition holds because   from Lemma \ref{L2.2} and Theorem \ref{T2.1},  calculation
for the canonical form of two arbitrary matrices can be reduced to the situation that at least one of the matrices is nonsingular. We further assume the second condition because otherwise  $v{\rm(P)}=-\infty$ (according to Remark \ref{R2.4}).

In general, numerical computation for Jordan normal form is sometimes unstable, see Chapter $7$ in \cite{golub2012matrix} and \cite{kaagstrom1980algorithm}. Large Jordan blocks are difficult to handle due to the unstableness in calculating Jordan blocks, see \cite{golub2012matrix,kaagstrom1980algorithm}.
 On the other hand, problem ${\rm(P)}$ itself also has some problematic cases, i.e., a small perturbation of the matrices will cause a significant change of the optimal solution, see \cite{pong2014generalized}, which, we believe, are strongly related to the unstableness of the Jordan decomposition methods. Our analysis above sheds some light on the relationship between the unstable cases of the Jordan normal form and the problematic cases of problem ${\rm(P)}$.
Fortunately, as we only need to calculate the Jordan normal form for real eigenvalues and  the blocks are of small sizes  of $1\times 1$ and $2\times 2$. Furthermore, if there are multiple $2\times2$ blocks, they must have the same eigenvalue. Thus, problem (P) is of a special structure in our implementation.
The steps in the  algorithms in appendices show that if the $2\times 2$ blocks are of the same eigenvalues and the number of the $2\times 2$ blocks are small, which covers all non-problematic cases of problem ${\rm(P)}$, the method to calculate the canonical form will be very fast, i.e., with time complexity $O(n^3)$.
In the literature, there are $O(n^3)$ Jordan decomposition methods, see \cite{beelen1990computational,kaagstrom1980algorithm}.

\section{Extension to equality constrained and interval bounded variants of GTRS}
 This section extends the usage of the canonical form and the SOCP reformulation to the  equality constrained  problem ${\rm (EP)}$ and the interval bounded problem  ${\rm (IP)}$.

\subsection{GTRS\ with equality constraint}
We first make an assumption similarly to Assumption \ref{A2.1}  to avoid some trivial cases.
\begin{asmp}\label{A3.1}
i) There is at least one feasible solution in ${\rm (EP)}$;
ii) The following conditions do not hold true at the same time: $A\succeq0$, $b\in {\rm Range}(A)$ and $c=\frac{1}{2}b^TA^{+}b$; and
iii) $A\neq0$.
\end{asmp}
This assumption is a ``two-side" Slater condition (Assumption 1 in \cite{xia2014s}) plus the condition $A\neq0$.
In fact, if $A=0,$ we can transform the constraint $b^Tx+c=0$  to another quadratic equality constraint which satisfies all the three conditions in ii), i.e., $(b^Tx+c)^2=0$.
So  problem ${\rm (EP)}$ can be transformed to an unconstrained quadratic optimization problem  using the null space representation of $L(x +A^+b)=0$, by decomposing $A$ as $A=L^TL,$   when  ii) is violated, or using the null space representation of  $b^Tx+c=0$ when iii) is violated.

With the same notations as in Section 2, Theorems \ref{T2.2}, \ref{T2.3} and \ref{T2.4} still hold here, which can be proved in a similar way by the S-lemma with equality \cite{xia2014s}.
However, Theorem \ref{T2.5} needs some modifications. In the following of this section, we still use $(\ref{cons})$ and $(\ref{obj})$ to denote the associated terms in both the  constraint
and objective functions, i.e., $(\ref{cons})=\frac{1}{2}\tau_1z^TE_Iz=\tau_1z_1z_2=\pi$ and $(\ref{obj})=\tau_1\lambda z_1z_2+\frac{1}{2}\tau_1z_2^2+e_1z_1+e_2z_2$.
\begin{thm}\label{T4.1}
Consider the case where there exists a type A block pair $(\tau_1E_1,\\ \tau_1E_1J_{1}(\lambda,2))$ in problem ${\rm (P)}$ and the eigenvalue of the associated Jordan
block $J_{1}(\lambda,2)$ is real. Assume there is a feasible solution $\bar{x}=(\bar{z}^T,\bar{y}^T)^T$ and let
$\pi=\tau_1\bar{z}_1\bar{z}_2$.
Let $\rho=\inf \{~(\ref{obj})  \mid(\ref{cons})=\tau_1z_1z_2=\pi\}$. We have the following three cases:
\begin{enumerate}
\item When $\tau_1=1$. If  $(e_1=0,e_2\neq0)$ or $(e_1=0,e_2=0, \pi=0)$, then $\rho=\lambda\pi-\frac{1}{2}e_2^2$ and the infimum is attainable;
\item  When $\tau_1=1$. If $e_1=0,e_2=0,\pi\neq0$, then $\rho=\lambda\pi-\frac{1}{2}e_2^2$ and the infimum is unattainable;
\item Otherwise, $\rho=-\infty$  and thus ${\rm(EP)}$ is unbounded from below.
\end{enumerate}
\end{thm}
\proof
The proof is similar to that of Theorem \ref{T2.5}.
\endproof

\begin{thm}\label{C3.1} If the optimal value of problem $\rm(EP)$ is bounded from below, then:
\begin{enumerate}
\item $\dim{E_i}\leq2, ~i=1,\ldots,p$, $\dim{E_i}=1,~i=p+1,\ldots,m$, and there is no complex eigenvalue pair in $J(\kappa_i,n_i)$;
\item If for some index $i$, $\dim{E_i}=2$, then the $i$th block satisfies case 1 or case 2 in Theorem \ref{T4.1}.
\end{enumerate}
\end{thm}

Note that the conditions in items 1 and 2 of Theorem \ref{C3.1} are necessary for problem (EP) to be
bounded from below and we assume that the conditions hold in the following of this section.  In the same way as the method in solving problem $\rm(P)$, we can then assume that $A$ and $D$ have the form in $(\ref{FinA})$ and $(\ref{FinB})$, $b_j=0$, for $j=l+1,\ldots,n$, and $e_{l+2j-1}=0,~j=1,\ldots,\frac{n-l}{2}$.

Similarly to Theorem \ref{mainthm},  using the S-lemma with equality \cite{xia2014s} under Assumption \ref{A3.1}, we have the following theorem.
\begin{thm}\label{T3.3}
Assume that items $1$ and $2$ in Theorem \ref{C3.1} are satisfied, then problem ${\rm(EP)}$ has the same optimal value with the following SOCP reformulation:
\begin{eqnarray*}
{\rm (EP_1)}~&~\min &\sum_{i=1}^l(\delta_i y_i+e_i x_i)+\sum_{j=1}^\frac{n-l}{2}\zeta_j z_j+c_0 \\
&~{\rm{s.t.}}&\sum_{i=1}^l(\alpha_i y_i+b_i x_i)+\sum_{j=1}^\frac{n-l}{2}z_j+c=0\\
&~&\frac{1}{2}x_i^2-y_i\leq0,~\forall i=1,2,\ldots,l,\\
&~&x,y\in \Re^l,~ z\in \Re^{\frac{n-l}{2}},
\end{eqnarray*}
where $c_0=-\sum_{j=1,\ldots,\frac{n-l}{2}}\frac{1}{2}e_{l+2j}^2$.

 More specifically, if  ${\rm(EP_1)}$ admits an optimal solution, then there exists an optimal solution $(\bar{x},\bar{y},\bar{z})$ to $\rm{(EP_1)}$
 with $\frac{1}{2}\bar{x}_i^2=\bar{y}_i,~i=1,2,\ldots,l$.
Moreover, we can find an optimal solution (or an $\epsilon$
optimal solution) $\tilde{x}$ to $(EP)$ with\begin{eqnarray*}
\left.\begin{array}{ll}
\tilde{x}_i=\bar{x}_i,~i=1,\ldots,l,\\
\tilde{x}_{l+2j}=\left\{\begin{array}{ll}
1/M&\text{if }e_{l+2j}=0,~\bar{z}_j\neq0 ,\label{une}\\
-e_{l+2j}&\text{otherwise},
\end{array}\right. ~~~~~~~~~~~~j=1,\ldots,\frac{n-l}{2},\\
\tilde{ x}_{l+2j-1}=\frac{\bar{z}_j}{\tilde{x}_{l+2j}},~j=1,\ldots,\frac{n-l}{2}.
\end{array}\right.
\end{eqnarray*}Particularly, if  ${\rm(EP)}$ is bounded from below, the optimal value of ${\rm(EP)}$ is unattainable if and only if $e_{l+2j}=0$ and $\bar{z}_j\neq0$.
In this case, for any $\epsilon>0$, there exists an  $\epsilon$ optimal solution $\tilde{x}$ such that $f(\tilde{x})-v{\rm(EP)}<\epsilon$ with a sufficient large $M>0$.

Furthermore, if $\zeta_i\neq \zeta_j$ for some $i\neq j$, where $i,j\in \{1,\ldots,\frac{n-l}{2}\}$,  problem $\rm{(EP)}$ is unbounded from below.
\end{thm}

The dual problem of ${\rm(EP_1)}$ has the same form with ${\rm(D_1)}$, except that $\nu\geq0$ is replaced by $\nu\in \Re$. From Theorem \ref{T3.3}, we can also conclude that all the $2\times2$ Jordan blocks
have the same eigenvalue, i.e., $\zeta_i=\zeta_j$ for $1\leq i,j\leq \frac{n-l}{2}$, if ${\rm(EP_1)}$ is bounded from below.
\begin{thm}\label{T3.4}
Under Assumption \ref{A3.1}, the objective values of $\rm{(EP_{1})}$ and the following Lagrangian dual
 problem are equal:\begin{eqnarray*}
({\rm{D_1}})~~
\max\limits_{\nu\in\Re}(c\nu+c_0+\sum_{i=1}^l h_i(\nu)+g(\nu)),
\end{eqnarray*}
where
\begin{eqnarray*}
h_i(\nu)=
\left\{\begin{array}{ll}
-\frac{(\nu b_i+e_i)^2}{2(\nu\alpha_i+\delta_i)}&\text{if }\nu\alpha_i+\delta_i>0,\\
0&\text{if }\nu\alpha_i+\delta_i=0\text{ and }e_i+\nu b_i=0,\\
-\infty&\text{otherwise},
\end{array}\right.
\end{eqnarray*}
$i=1,\ldots,l,$ and
\begin{eqnarray*}
g(\nu)=\left\{\begin{array}{ll}
0&\text{if }\ \zeta_j+\nu=0,~\forall j=1,\ldots,\frac{n-l}{2},\\
-\infty&\text{otherwise}.
\end{array}\right.
\end{eqnarray*}
\end{thm}

Assumption \ref{A3.1} is necessary in the above theorem, since otherwise ${(\rm EP_1)}$ may not be equivalent to ${(\rm EP)}$.
In fact, $v{(\rm EP_1)}=-\infty$ and $v{(\rm EP)}>-\infty$ hold true for the following case stated in Theorem 3 in \cite{xia2014s}:
$D$ has exactly one negative eigenvalue, $A=0$, $b\neq0$ and
$$\left(\begin{array}{ccc}
V^TDV&  V^T(Dx_0+e)\\
(x_0^TD+e^T)V &f(x_0)
\end{array}
\right)\succeq0,$$
where $x_0=-\frac{e}{2b^Tb}b$, $V\in \Re^{n\times(n-1)}$ is the matrix basis of the null space  of the $1\times n$ matrix $b^T$.
For example, consider the problem $\min \{f(x)\mid h(x)=0\}$, where  $f(x)=2x_1^2-x_2^2$ and $h(x)=x_1-x_2$.
In this case $v{(\rm EP_1)}=-\infty$ but $v{(\rm EP)}=0$.

Now let us denote $h(x)=\frac{1}{2}x^TAx+b^Tx+c$ and $f(x)=\frac{1}{2}x^TDx+e^Tx+v$ in the following of this section. Recently, Xia et al.
\cite{xia2014s} give the conditions when the S-lemma with equality holds, with or without Assumption 1 in \cite{xia2014s} (the ``two-side" Slater condition).
The S-lemma with equality asks if the following two statements are equivalent:\\
$~~~~{\rm(E_1)}$ $(\forall x\in \Re^n)$ $h(x)=0\Rightarrow f(x)\geq0$.\\
$~~~~{\rm(E_2)}$ $\exists$ $\mu\in \Re$ such that $f(x)+\mu h(x)\geq0, ~\forall x\in \Re^n$. \\
Theorem 3 in \cite{xia2014s} states that under the ``two-side" Slater condition, ${\rm(E_1)}$ is equivalent to ${\rm(E_2)}$ except for the following case:
$D$ has exactly one negative eigenvalue, $A=0$, $b\neq0$ and
$$\left(\begin{array}{ccc}
V^TDV&  V^T(Dx_0+e)\\
(x_0^TD+e^T)V &f(x_0)
\end{array}
\right)\succeq0,$$
where $x_0=-\frac{c}{2b^Tb}b$ and $V\in \Re^{n\times(n-1)}$ is the matrix basis of the null space of the $1\times n$ matrix $b^T$.
In fact, the above case violates condition iii) in Assumption \ref{A3.1}. So under our assumption, this case is out of consideration and the S-lemma with equality always holds in our investigation.

Under Assumption \ref{A3.1}, by applying the canonical form and the SOCP reformulation of ${\rm(EP)}$, we obtain a simplified S-lemma with equality, which shows the following two statements are equivalent:\\
$~~~~{\rm(E_1)}$  ($\forall x\in \Re^n$) $h(x)=0\Rightarrow f(x)\geq0$.\\
$~~~~{\rm(\tilde{E}_2)}$ $\exists x,y\in \Re^l, ~z\in \Re^{\frac{n-l}{2}}:\left\{\begin{array}{ll}
\delta^T y+e^T x+\zeta^T z+c_0+v\geq0\\
\alpha^Ty+b^Tx+\bar{1}^Tz+c=0\\
\frac{1}{2}x_i^2-y_i\leq0,~\forall i=1,\ldots,l.\end{array}\right.$\\
The equivalence between ${\rm(E_1)}$  and ${\rm(\tilde{E}_2)}$ can be directly derived from Theorem \ref{T3.3}.
One advantage of  our simplified S-lemma with equality is its feature of LMI free, thus being  more tractable\ for large-scale problems, while the S-lemma in \cite{xia2014s} is not.

\subsection{GTRS with interval constraint}
Similarly to the  equality constrained case, we make the following assumption.
\begin{asmp}\label{A3.2}
i) There is at least one feasible solution in ${\rm (IP)}$;
ii) The following conditions do not hold true at the same time: $A\succeq0$, $b\in {\rm Range}(A)$ and $c_i=\frac{1}{2}b^TA^{+}b$ for $i=1$ or $2$; and
iii) $A\neq0$.
\end{asmp}

Theorem \ref{C3.1} holds in this case and thus we can still assume, without loss of generality, $A$ and $D$ have the form in $(\ref{FinA})$
and $(\ref{FinB})$, $b_j=0$, for $j=l+1,\ldots,n$, and $e_{l+2j-1}=0,~j=1,~\ldots,~\frac{n-l}{2}.$
\begin{thm}

\label{T4.3}
Assume that the conditions in items $1$ and $2$ in Theorem \ref{C3.1} are satisfied, problem ${\rm(IP)}$ has the same optimal value with the following SOCP problem:
\begin{eqnarray*}
{\rm (IP_1)}~&~\min &\sum_{i=1}^l(\delta_i y_i+e_i x_i)+\sum_{j=1}^\frac{n-l}{2}\zeta_j z_j+c_0 \\
&~{\rm{s.t.}}&\bar{h}(x,y,z)=\sum_{i=1}^l(\alpha_i y_i+b_i x_i)+\sum_{j=1}^\frac{n-l}{2}z_j\geq c_1\\
&~&\bar{h}(x,y,z)=\sum_{i=1}^l(\alpha_i y_i+b_i x_i)+\sum_{j=1}^\frac{n-l}{2}z_j\leq c_2\\
&~&\frac{1}{2}x_i^2-y_i\leq0,~\forall i=1,2,\ldots,l,\\
&~&x,y\in \Re^l, ~z\in \Re^{\frac{n-l}{2}},
\end{eqnarray*}
$c_0=-\sum_{j=1,\ldots,\frac{n-l}{2}}\frac{1}{2}e_{l+2j}^2$.

 More specifically, if  ${\rm(IP_1)}$ admits an optimal solution, then there exists an optimal solution $(\bar{x},\bar{y},\bar{z})$ to $\rm{(IP_1)}$
 with $\frac{1}{2}\bar{x}_i^2=\bar{y}_i,~i=1,2,\ldots,l$.
Moreover, we can find an optimal solution (or an $\epsilon$
optimal solution) to (IP) with
\begin{eqnarray*}
\left.\begin{array}{ll}
\tilde{x}_i=\bar{x}_i,~i=1,\ldots,l,\\
\tilde{x}_{l+2j}=\left\{\begin{array}{ll}
1/M&\text{if }e_{l+2j}=0,~\bar{z}_j\neq0 ,\\
-e_{l+2j}&\text{otherwise},
\end{array}\right. ~~~~~~~~~~~~j=1,\ldots,\frac{n-l}{2},\\
\tilde{ x}_{l+2j-1}=\frac{\bar{z}_j}{\tilde{x}_{l+2j}},~j=1,\ldots,\frac{n-l}{2}.
\end{array}\right.
\end{eqnarray*}
Particularly, if  ${\rm(IP)}$ is bounded from below, the optimal value of ${\rm(IP)}$ is unattainable if and only if $e_{l+2j}=0$ and $\bar{z}_j\neq0$.
In this case, for any $\epsilon>0$, there exists an  $\epsilon$ optimal solution $\tilde{x}$ such that $f(\tilde{x})-v{\rm(IP)}<\epsilon$ with a sufficient large $M>0$.

Furthermore, if $\zeta_i \neq \zeta_j$ for some $i\neq j$, where $i, j \in \{1,\ldots,\frac{n-l}{2}\}$,  problem $\rm{(IP)}$ is unbounded from below.
\end{thm}

\proof
${\rm(IP)}$ is equivalent to the following ${\rm(IP_2)}$:
\begin{eqnarray*}
{\rm (IP_2)}~&~\min &\sum_{i=1}^l(\delta_i y_i+e_i x_i)+\sum_{j=1}^\frac{n-l}{2}\zeta_j z_j+c_0 \\
&~{\rm{s.t.}}&-\sum_{i=1}^l(\alpha_i y_i+b_i x_i)-\sum_{j=1}^\frac{n-l}{2}z_j+c_1\leq0\\
&~&\sum_{i=1}^l(\alpha_i y_i+b_i x_i)+\sum_{j=1}^\frac{n-l}{2}z_j- c_2\leq0\\\
&~&\frac{1}{2}x_i^2-y_i=0,~\forall i=1,2,\ldots,l,\\
&~&x,y\in \Re^l, ~z\in \Re^{\frac{n-l}{2}}.
\end{eqnarray*}
So we only need to prove the equivalence between $\rm(IP_1)$ and $\rm(IP_2)$.

By the S-lemma with interval bounds \cite{wang2014strong}, similarly to Lemma \ref{L2.3}, we know that if $v\rm(IP_1)$ is unbounded from below, then $v\rm(IP_2)$ is unbounded from below.

Now we consider the case where $v\rm(IP_1)$ is bounded from below. Then there exists a global minimum $(x^{*}, y^{*},z^{*})$ for  $\rm(IP_1)$.
The Fritz-John conditions of $\rm(IP_1)$ are stated as following:
There exist $\nu_0\geq0,~\nu_1\geq0,~\nu_2\geq0,~\mu_i\geq0,  ~i=1,\ldots,l$, not all of which are zero, such that
\begin{eqnarray*}
&~&~\nu_0\delta_i-(\nu_1-\nu_2)\alpha_i-\mu_i=0, ~\forall i=1,\ldots,l,\\
&~&~\nu_0e_i-(\nu_1-\nu_2)b_i+\mu_ix_i^{*}=0, ~\forall i=1,\ldots,l,\\
&~&~\nu_0\zeta_j-(\nu_1-\nu_2)=0, ~\forall j=1,\ldots,\frac{n-l}{2}.
\end{eqnarray*}
We assume that $\alpha_i$ and $b_i$ are not both zero for $i=1,\ldots,l$, otherwise $x_i$ is a free variable only appearing in the objective function and then $\rm(IP_{1})$ is either
unbounded from below or can be reduced to a new problem without variable $x_i$. Moreover, we cannot take equality in both sides of the quadratic
 constraint, so there must exist at least one strict inequality. Then, from the last equation in Fritz-John
 conditions and the  complementary slack conditions, we conclude $\nu_1(\bar{h}(x^{*},y^{*},z^{*})-c_1)=0,~\nu_2(\bar{h}(x^{*},y^{*},z^{*})-c_2)=0,$  and one of
  $\nu_1$ and $\nu_2$ must be $0$.  Then from the first equation in Fritz-John conditions we know that if there exists some index $i$ such that
   $\frac{1}{2}(x_i^*)^2-y_i^{*}<0$, together with the  complementary slack conditions  $\mu_i(\frac{1}{2}(x_i^*)^2-y_i^{*})=0$, we conclude $\mu_i=0$ and $\nu_0>0$ (otherwise  $(\nu_1-\nu_2)\alpha_i=(\nu_1-\nu_2)b_i=0 \Rightarrow \nu_1-\nu_2=0 \Rightarrow \nu_1=\nu_2=0$ and thus $\mu_i=0$ for all $i$, which contradicts the fact
   that $\nu_0$, $\nu_1$, $\nu_2$, $\mu_i$,  $i=1\cdots,l$ are not all zero).
 So the Fritz-John conditions is reduced to the KKT conditions. Because one of $\nu_1$ and $\nu_2$ is $0$,  Assumption 6 in \cite{ben2014hidden} holds.
 Then, by applying Theorem 7 in \cite{ben2014hidden}, we can get another optimal solution $(\bar{x},\bar{y},\bar{z})$ to $\rm{(IP_{1})}$
 with $\frac{1}{2}\bar{x}_i^2=\bar{y}_i,~i=1,2,\ldots,l$, and $\bar{z}=z^*$, which is also optimal to $\rm{(IP_{2})}$.

 The remaining of the proof is similar to that of Theorem \ref{mainthm}.
\endproof

Similarly, the dual problem of $\rm(IP_1)$ is a simple problem with two variables.
\begin{thm}
Under Assumption \ref{A3.1}, the objective values of $(\rm IP_1)$ and the following Lagrangian dual
problem are equal:\begin{eqnarray*}
({\rm ID_1})
\max\limits_{\nu\geq0}(\nu_1c_1-\nu_2c_2+c_0+\sum_{i=1}^lh_i(\nu)+g(\nu)),
\end{eqnarray*}
where
\begin{eqnarray}
h_i(\nu)&=&
\left\{\begin{array}{ll}
-\frac{(-\nu_1b_i+{\nu}_2b_i+e_i)^2}{2(-{\nu}_1\alpha_i+{\nu}_2\alpha_i+\delta_i)}&\text{if }{-\nu}_1\alpha_i+{\nu}_2\alpha_i+\delta_i>0,\\
0&\text{if }{-\nu}_1\alpha_i+{\nu}_2\alpha_i+\delta_i=0\ \text{and}\ e_i-{\nu}_1b_i+{\nu}_2b_i=0,\\
-\infty&\text{otherwise},
\end{array}\right. \nonumber\\
&~&~~~~~~~~~~~~~~~~~~~~~~~~\label{h_1-2}
\end{eqnarray}
$i=1,\ldots,l,$ and
\begin{eqnarray}
g(\nu)=\left\{\begin{array}{ll}
0&\text{if}\ \zeta_j-\nu_1+\nu_2=0,\text{ for all}\,j=1,\ldots,\frac{n-l}{2},\\
-\infty&\text{otherwise}.
\end{array}\right. \label{g-2}
\end{eqnarray}
\end{thm}

\proof
The Lagrangian of $\rm(IP_1)$ is:
\begin{eqnarray*}
L(x,y,z;\mu,\nu)&=&\delta^T y+e^T x+\zeta^T z+c_0+\nu_1(-\alpha^Ty-b^Tx-\bar{1}^Tz+c_1)\\
&~&+\nu_2(\alpha^Ty+b^Tx+\bar{1}^Tz-c_2)+\sum_{i=1}^l\mu_i(\frac{1}{2}x_i^2-y_i)\\
&=&\sum_{i=1}^ly_i(\delta_i+\nu_2\alpha_i-\nu_1\alpha_i-\mu_i)+\sum_{i=1}^lx_i(e_i+\frac{1}{2}\mu_ix_i+\nu_2b_i\\
&~&-\nu_1b_i)+\sum_{j=1}^{\frac{n-l}{2}}(\zeta_j+\nu_2-\nu_1) z_j+c_1\nu_1-c_2\nu_2+c_0.
\end{eqnarray*}

The dual Lagrangian function is then read as
\begin{eqnarray*}
m(\mu,\nu)=\min\limits_{x,y,z}L(x,y,z;\mu,\nu).
\end{eqnarray*}
As the objective function is separable, we can get an analytical form of $m(\mu,\nu)$:
\begin{eqnarray*}
\max\limits_{\mu\geq0,\nu\geq0}m(\mu,\nu)=\max\limits_{\nu\geq0}(c_1\nu_1-c_2\nu_2+c_0+\sum_{i=1}^lh_i(\nu)+g(\nu)),
\end{eqnarray*}
where $h_i(\nu)$ and $g(\nu)$ are given in (\ref{h_1-2}) and (\ref{g-2}), respectively.
\endproof

Since one of the two Lagrangian multipliers $\nu_1$ and $\nu_2$ must be zero, we can separate problem $\rm(ID_1)$ into two problems with a single variable $\nu_1$ or $\nu_2$ by setting either $\nu_1=0$ or $\nu_2=0$. Then problem   $\rm(ID_1)$ is reduced to a problem with the same form as problem ${\rm(D_1)}$.
Besides,  problem $\rm(ID_1)$ is also equivalent to an SOCP problem with two variables similar as problem ${\rm(D_1)}$.

\begin{rem}
Actually, problem ${\rm(IP)}$ must have an optimal solution on the boundary, except for the case where $D\succeq0$, $e\in Range(D)$ and $x=-D^+e$ is in the interior
of the interval constraint. This is because, if the optimal solution $x^*$ is not on the boundary, then  $x^*$ must be a local minimum and  ${\rm(IP)}$  has only one local minimum under
the conditions that
$D$ is semi-definite positive, $e\in Range(D)$ and the local minimum is $x=-D^+e$. So we can first verify whether the conditions $D\succeq0$ and $e\in Range(D)$ are satisfied and then check whether $x=-D^+e$ is in the interior
of the constraint. Otherwise, the optimal solution must be on the boundary.  Then we can separate the problem into two equality constrained problems with an equality constraint
$h(x)=c_1$ or $h(x)=c_2$, and solve them  with the  methods for the equality constrained case.
And the solution with the smaller optimal value of the above two equality constrained problems is the optimal solution of problem ${\rm(IP)}$. \end{rem}

We now denote $h(x)=\frac{1}{2}x^TAx+b^Tx$ and $f(x)=\frac{1}{2}x^TDx+e^Tx+v$. Similarly to the simplified S-lemma with equality, we can simplify the S-lemma with interval bounds in \cite{wang2014strong}
under Assumption \ref{A3.2} to the following form: ${\rm(I_1)}$ $\Leftrightarrow$ ${\rm(\tilde{I}_2)}$. \\
$~~~~{\rm(I_1)}$  ($\forall x\in \Re^n$) $c_1\leq h(x)\leq c_2\Rightarrow f(x)\geq0$.\\
$~~~~{\rm(\tilde{I}_2)}~\exists x,y\in \Re^l, ~z\in \Re^{\frac{n-l}{2}}:\left\{\begin{array}{ll}
\delta^T y+e^T x+\zeta^T z+c_0+v\geq0,\\
c_1\leq\alpha^Ty+b^Tx+\bar{1}^Tz\leq c_2,\\
\frac{1}{2}x_i^2-y_i\leq0,~\forall i=1,\ldots,l.\end{array}\right.$
\\

\section{Conclusions}
In this study, we have successfully developed an SOCP reformulation for the GTRS, which is in fact a quadratic programming over a single nonconvex quadratic constraint. Particularly, we have derived the SOCP\ reformulation under the condition that the GTRS  is bounded from below, via a canonical congruent form of the two matrices in both the objective and constraint
 functions. While Ben-Tal and Hertog investigate in \cite{ben2014hidden} the simultaneous diagonalizability of the two matrices, we explore the simultaneous block diagonalizability,
 which applies to arbitrary two matrices.
 More specifically, we introduce and extend the canonical form for
  two real matrices in \cite{uhlig1976canonical} to find a block diagonal form for two matrices  in both the objective and constraint functions. Exploiting the
  separability of the block diagonal form of the two matrices, we show that problem (P) is SOCP representable if it is bounded from below.
  We also establish the attainableness of the problem from the canonical form without any additional calculation.   Moreover, as our SOCP reformulation is LMI free, it can be solved much faster than the SDP reformulation for the GTRS. We further extend our methods to solve  variants of problem ${\rm(P)}$
where the inequality constraint is replaced by either an equality constraint or an interval constraint.
Moreover, we obtain simplified versions of  the S-lemma under the three kinds of constraints.

One of our future research is to consider variants of the GTRS with additional linear inequality constraints or two general quadratic constraints.

\begin{appendices}
In the appendices, we present an algorithm to compute the canonical form if all the eigenvalues of the Jordan matrix of $A^{-1}D$ ($A$ is assumed invertible)
are real and there only exist $1\times 1$ and $2\times 2$ Jordan blocks, where all the $2\times 2$ Jordan blocks have the same eigenvalue.
Without loss of generality, we assume the non-singularity of matrix $A$ for simplicity of the analysis.
\begin{dfn}
A single element $a$, a $2\times1 $ matrix $\left(\begin{array}{c}0\\a\end{array}\right)$, a $1\times 2 $ matrix $(0,a)$  or a $2\times 2$ matrix $\left(\begin{array}{cc}0&a\\a&b\end{array}\right)$ is called a lower striped matrix, where $a,b \in \Re$.
\end{dfn}
\begin{dfn}
Let $J(\lambda,n_1)$, $J(\lambda,n_2)$, $\dots$, $J(\lambda,n_l)$ denote all the Jordan blocks associated with the same eigenvalue $\lambda$ of a real matrix A. Then $C(\lambda)={\rm diag}(J(\lambda,n_1),\ldots,J(\lambda,n_l))$, where ${\rm dim}J_i\geq{\rm dim}J_{i+1}$ for $i$ = 1, $\ldots$ $l-1$, is called the full chain of Jordan blocks or full Jordan chain of length l associated with $\lambda$.
\end{dfn}

We present the following two algorithms to calculate the canonical form. Algorithm 1 is the main algorithm to compute the canonical form. Algorithm 2 is a subroutine to calculate the canonical form for the full Jordan chain associated with $\lambda$.

Let us look into Algorithm 1 first. {Line $2$ is to find the Jordan normal form of $A^{-1}D$, i.e., $J=V^{-1}A^{-1}DV$, where $V$ is some invertible matrix. Lines 3--5 are just to avoid some unbounded cases of problem ${\rm(P)}$.}
Line $6$ updates $A$ and $D$: $A=V^{T}AV={\rm diag}(A_1,A_2,\cdots,A_m)$ and $D=V^{T}DV={\rm diag}(D_1,D_2,\cdots,D_m)$.
According to the proof of Theorem 1 in \cite{uhlig1976canonical}, $A$ and $D$ are now both block diagonal matrices and the blocks $A_i$ and $D_i$ have the same dimension with $C(\lambda_i)$, (in fact $A_i^{-1}D_i=C(\lambda_i)$), where $C(\lambda_i)$ is the full Jordan chain associated with $\lambda_i$, and has the form $C(\lambda_i)={\rm diag}(J_1(\lambda_i),J_2(\lambda_i),\cdots,J_{l_i}(\lambda_i))$, where $$J_1(\lambda_i)=\ldots=J_k(\lambda_i)=\left(\begin{array}{cc}\lambda_i&1\\0&\lambda_i\end{array}\right)\text{ and }J_{k+1}(\lambda_i)=\ldots=J_{l_{i}}(\lambda_i)=\lambda_i.$$
Note that from the assumption there is at most one Jordan chain that has $2\times 2$ Jordan blocks and we further assume that this Jordan chain is associated with the eigenvalue $\lambda_i$. According to the Lemma 2 in \cite{uhlig1976canonical}, $A_i$ is a block matrix in the following form,
 $$\left(\begin{array}{cccc}
A_{11}^i&\cdots&A_{1l}^i\\
\vdots&&\vdots\\
A_{l1}^i&\cdots&A_{ll}^i
\end{array}\right),$$ where $A_{pq}^i$ is a $\dim J_p(\lambda_i)\times \dim J_q(\lambda_i)$ lower striped matrix, $p,q=1,\ldots,l$ (we use $l=l_i$ here and in the following for simplicity).
\begin{algorithm}[!ht]
  \caption{Calculate the canonical form for two real symmetric matrices}
  \begin{algorithmic}[1]
  \Require $A$ and $D$ are $n\times n$ real symmetric matrices and $A$ is nonsingular  \Ensure $\check{A}$ and $\check{D }$ in the canonical form and the congruent matrix $U$ and the final matrices $\check{A}= U^TAU$ and $\check{ D}= U^TDU$
  \Function {Canonical-transformation}{$A, D$}
    \State Calculate the Jordan normal form for $A^{-1}D$, find a nonsingular matrix $V$ such that $J=V^{-1}A^{-1}DV={\rm diag}(C(\lambda_1),\ldots,C(\lambda_m))$, where $C(\lambda_i)$ is the full
     Jordan chain associated with the eigenvalue $\lambda_i$, denote $l_i$  the number of blocks on the diagonal of $C(\lambda_i)$ and  $k_i$ the number of $2\times2 $
     blocks.\footnotemark
    \If {$J$ has complex eigenvalues or Jordan blocks with size greater than $2$ or more than one full Jordan chain has $2\times 2$ Jordan blocks}
    \State \Return{(in this case problem ${\rm(P)}$ is unbounded from below)}
    \EndIf
 \State $A=V^{T}AV={\rm diag}(A_1,A_2,\cdots,A_m)$ and $D=V^{T}DV={\rm diag}(D_1,D_2,\cdots,D_m)$
    \For{$i=1$; $i\leq m$; $i++$ }
        \State $(\check{A}_i,\check{D}_i,U_i)=$\Call{JBF}{$A_i,D_i,C(\lambda_i),k_{i},l_{i}$}
    \EndFor
    \State $\check{A}= {\rm diag}(\check{A}_1,\ldots,\check{A}_m)$, $\check{D}= {\rm diag}(\check{D}_1,\ldots,\check{D}_m)$, $U= {\rm diag}(U_1,\ldots,U_m)$
    \EndFunction
      \end{algorithmic}

\end{algorithm}

\begin{algorithm}[!ht]
  \caption{Calculate the canonical form for a Jordan chain }
  \begin{algorithmic}[1]
  \Require $A$ and $D$ are $l\times l$ real symmetric matrices and $A$ is nonsingular, $J$ is the Jordan matrix of $A^{-1}D$ and has only $1\times 1$ and $2\times 2$ Jordan blocks with a real eigenvalue, where the first $k$ blocks are $2\times2 $ blocks ($k$ could be equal to $0$) and the last $l-k$ blocks are $1\times1$ blocks  \Ensure $\check{A}$ and $\check{D}$ in the canonical form and the congruent matrix $U$ such that $\check{ A}=U^TAU$ and $\check{D}=U^TDU$
    \Function {$(\check{A},\check{D},U)$$=$ JBF}{$A,D,J,k,l$}~~~~
          \State $U= I_{l+k}$
      \For {$j=1$; $j\leq k$; $j++$ }
         \State $X_j= I_{l+k}$
         \If {the $(j,j)$th block of $A$, $A_{jj}$ is singular}
         \State $X_j=(\ref{X1})$
         \Else
         \State$X_j=(\ref{X2})$
         \EndIf
         \State $Y_{j}=(\ref{Y}),~U= UX_{j}Y_j $
      \EndFor
      \State $\bar{A}= U^TAU$
      \State Find the associated congruent matrix $Q$ of the spectral decomposition of lower right  $(l-k)\times(l-k)$ submatrix $\bar{A}_{(k+1):l,(k+1):l}$
      \State $\bar{Q}={\rm diag}(I_{2k},Q)$, $\tilde{A}= \bar{Q}^T\bar{A}\bar{Q}$
      \For {$j=1;j\leq k;j++$}
      \State $P_j=\left(\begin{array}{cc}
\frac{1}{\sqrt{|a_{j1}|}}&-\frac{a_{j2}}{2a_{j1}\sqrt{|a_{j1}|}}\\
0&\frac{1}{\sqrt{|a_{j1}|}}
\end{array}\right)(\text{Note the }j\text{th block of } \tilde{A} \text{ is} \left(\begin{array}{cc}0&a_{j1}\\a_{j1}&a_{j2}\end{array}\right))$
      \EndFor
      \State $P_{k+1}={\rm diag}(\frac{1}{\sqrt{|a_{k+1}|}},\cdots,\frac{1}{\sqrt{|a_l|}})$ (where $a_i$ is the $(2k+i)$th element in the diagonal of $\tilde{A}$)
      \State $P={\rm diag}(P_1,\ldots,P_k,P_{k+1}), ~\check{A}=P^T\tilde{A}P,\check{D}=\check{A}J,~U=U\bar{Q}P$
    \EndFunction
  \end{algorithmic}
\end{algorithm}

On Lines 7--10 of Algorithm 1,
we switch to Algorithm 2 for the Jordan chain $C(\lambda_i)$ and then on Line 10 we
use diagonal matrix $U= {\rm diag}(U_1,\ldots,U_m)$ as a congruent matrix to compute the canonical form of $A$ and $B$.

Next we describe Lines 2--12 in Algorithm 2.
Assuming the first $k$ Jordan blocks in  $C(\lambda_i)$ are of size $2\times 2$ (we use $k=k_i$ here and in the following for simplicity), we next make a transformation to make $A_{11}^i$ nonsingular if $A_{11}^i$ is singular. If $A_{11}^i$ is singular, we show that  after a suitable permutation to $A_{i}$,
\footnotetext{In the algorithm, we arrange all  the  $2\times2 $ blocks to the upper left part of $C(\lambda_i)$ and all $1\times1$ blocks to
the lower right part when calculating the Jordan normal form.}
the new $A_{11}^i$ will be nonsingular.
If there exists  some $s,~2\leq s\leq k$ such that $A_{ss}^i$ is nonsingular, let
\begin{eqnarray*}
X_i^1={\rm diag}(\left(\begin{array}{cccccccc}
0&0&\cdots&0&I_2\\
0&I_2&&&\\
\vdots&&\ddots&&\\
0&&&I_2&\\
I_2&&&&0
\end{array}\right)_{2s\times 2s},~I_2,\ldots,I_2,~I_{l-k})\end{eqnarray*}
then $(X_i^1)^TA_iX_i^1$ has a nonsingular $(1,1)$ block. Otherwise all $A_{ss}^i$ are singular,
$s=1,\ldots,k$. Since $A_i$ is nonsingular and $A_{11}^i, A_{12}^i,\ldots,A_{1l}^i$ are low striped matrices, there must be at least one nonsingular matrix among
$ A_{12}^i,\ldots,A_{1k}^i$ (otherwise the first row of $A_i$ is $0$, which contradicts the nonsingularity of $A_i$).
Suppose that $A_{1s}^i$ is nonsingular for $2\leq s\leq k$ and let
\begin{eqnarray*}
X_i^1=I_{l+k}+\left(\begin{array}{cccc}
0_{2\times 2}&0_{2\times(2s-2)}&-I_2&0_{2\times(l+k-2s)}\\
0_{(2s-2)\times2}&\ddots&0_{(2s-2)\times2}&\\
I_2&0_{2\times(2s-2)}&0_{2\times2}&\\
0_{(l+k-2s)\times2}&&&
\end{array}\right).
\end{eqnarray*}
Then $$(X_i^1)^TA_iX_i^1=\left(\begin{array}{cc}
A_{11}^i+A_{1s}^i+(A_{1s}^i)^{T}+A_{ss}^i&*\\
*&*\end{array}\right)$$ with nonsingular $(1,1)$ block.
Now, as $A_i=(X_i^1)^TA_i{X_i^1}$ has nonsingular $(1,1)$ block, we can transform $A_i$ next to a simpler matrix with only $0$ blocks in $(1,t)$ and $(t,1)$ positions via congruence, for all  $t = 2, \ldots, l$.
Let
\begin{eqnarray*}
Y_i^1=I_{l+k}+\left(\begin{array}{cccc}
0&-(A_{11}^i)^{-1}A_{12}^i&\cdots&-(A_{11}^i)^{-1}A_{1l}^i\\
0&0&\cdots&0\\
\vdots&\vdots&&\vdots\\
0&0&\cdots&0
\end{array}\right).
\end{eqnarray*}
Then  we get $$A_i^1=(Y_i^1)^TA_iY_i^1=\left(\begin{array}{cccc}
A_{11}^{i,1}&0&\cdots&0\\
0&A_{22}^{i,1}&\cdots&A_{2l}^{i,1}\\
0&\vdots&&\vdots\\
0&A_{l2}^{i,1}&\cdots&A_{ll}^{i,1}
\end{array}\right).$$

Next, we conduct similar operations to the sub-matrix,$$\left(\begin{array}{ccc}
A_{22}^{i,1}&\cdots&A_{2l}^{i,1}\\
\vdots&&\vdots\\
A_{l2}^{i,1}&\cdots&A_{ll}^{i,1}
\end{array}\right)$$
and then conduct similar operations again iteratively
for $k-2$ times with the associated congruent matrices in the $j$th iteration. Let us be more specific using similar notations as in the first iteration. If some $(s_j,s_j)$ block is nonsingular, $j+1\leq s_j\leq k$, we introduce
\begin{eqnarray}
X_i^j={\rm diag}(\underbrace{{I_2,\ldots,I_2}}_{j-1},K,I_2,\ldots,I_2,I_{l-k}),\label{X1}\end{eqnarray}
where $$K=\left(\begin{array}{cccccccc}
0&0&\cdots&0&I_2\\
0&I_2&&&\\
\vdots&&\ddots&&\\
0&&&I_2&\\
I_2&&&&0
\end{array}\right)_{2s_{j}\times 2s_{j}}.$$
Otherwise,
we denote \begin{eqnarray}
X_i^{j}=I_{l+k}+\left(\begin{array}{ccccc}
0_{(2j-2)\times(2j-2)}&0&0&0&0\\
0&0_{2\times 2}&0_{2\times(2s_{j}-2)}&-I_2&0_{2\times m}\\
0&0_{(2s_{j}-2)\times2}&\ddots&0_{(2s_{j}-2)\times2}&0\\
0&I_2&0_{2\times(2s_{j}-2)}&0_{2\times2}&0\\
0&0_{m\times2}&0&0&0
\end{array}\right),\label{X2}
\end{eqnarray}
where $m=l+k-2s_{j}-2j$, and \begin{eqnarray}
Y_i^j=I_{l+k}+\left(\begin{array}{ccccc}
0_{(2j-2)\times(2j-2)}&&&&\\
&0_{2\times2}&-(A_{22}^i)^{-1}A_{23}^i&\cdots&-(A_{22}^i)^{-1}A_{2l}^i\\
&0&0&\cdots&0\\
&\vdots&\vdots&&\vdots\\
&0&0&\ldots&0
\end{array}\right).\label{Y}
\end{eqnarray}
 Then we get$\bar{A}_i=A_i^k=(X_i^1Y_i^1\ldots X_i^kY_i^k)^TA_i(X_i^1Y_i^1\ldots X_i^kY_i^k)=$\\
$$ \left(\begin{array}{ccccccc}
A_{11}^{i,k}&&&&&\\
&\ddots&&&&\\
&&A_{kk}^{i,k}&&&\\
&&&A_{k+1,k+1}^{i,k}&\ldots&A_{k+1,l}^{i,k}\\
&&&\vdots&&\vdots\\
&&&A_{l,k+1}^{i,k}&\ldots&A_{l,l}^{i,k}
\end{array}\right),$$
where $A_{jj}^{i,k}$ is a $2\times 2$ low striped matrix for $j=1,\ldots,k$, and $A_{pq}^{i,k}$ is a single number for $p,q>k+1$.
(We rewrite  in the following $A_{jj}^{k}=A_{jj}^{i,k}$ for simplicity.)

Lines $13-14$ of Algorithm 2 state that: using spectral decomposition to the following submatrix, $$\bar{A}_{(k+1):l,(k+1):l}=\left(\begin{array}{ccc}A_{k+1,k+1}^k&\ldots&A_{k+1,l}^k\\
\vdots&&\vdots\\
A_{l,k+1}^k&\ldots&A_{l,l}^k
\end{array}\right),$$ we finally get $Q_i^T\bar{A}_{(k+1):l,(k+1):l}Q_i$=${\rm diag}(a_{k+1},\ldots,a_l)$. Denote $$W_i=X_i^1Y_i^1\ldots X_i^kY_i^k \text{ and } \bar{Q}_i={\rm diag}(I_{2k},Q_i),$$ and we get
$$\tilde{A}_i=\bar{Q}_i^T\bar{A}_i\bar{Q}_i={\rm diag}(A_{11}^{k},\ldots,A_{kk}^{k},a_{k+1},\ldots,a_l).$$
Since both $W_i$ and $\bar{Q}_i$ commute with $C(\lambda_i)$ according to \cite{uhlig1976canonical},
\begin{eqnarray*}
\tilde{D}_i&=&(W_i\bar{Q}_i)^TD_iW_i\bar{Q}_i\\
&=&(W_i\bar{Q}_i)^TA_iC(\lambda_i)W_i\bar{Q}_i\\
&=&(W_i\bar{Q}_i)^TA_iW_i\bar{Q}_iC(\lambda_i)\\
&=&\tilde{A}_iC(\lambda_i)\\
&=&{\rm diag}(D_{11}^{k},\ldots,D_{kk}^{k},d_{k+1},\ldots,d_l)
\end{eqnarray*}
has the same block diagonal form with $\tilde{A}$.

The remaining of Algorithm 2 considers the transformation of the following block to the canonical form,
$$A_{jj}^{k}=\left(\begin{array}{cc}
0&a_{j1}\\
a_{j1}&a_{j2}
\end{array}\right).$$
Introducing the following congruent  matrix
$$P_i^j=\left(\begin{array}{cc}
\frac{1}{\sqrt{|a_{j1}|}}&-\frac{a_{j2}}{2a_{j1}\sqrt{|a_{j1}|}}\\
0&\frac{1}{\sqrt{|a_{j1}|}}
\end{array}\right),$$ we get
$$(P_i^j)^TA_{jj}^kP_i^j=\epsilon_i^j E_i^j,$$
where  $E_i^j=\left(\begin{array}{cc}
0&1\\
1&0
\end{array}\right)$, $\epsilon_i^j=1$ if $a_{j1}>0$ or otherwise $\epsilon_i^j=-1$. Then from
$(P_i^j)^TA_{jj}^kP_i^j=\epsilon_i^j E_i^j$, we know
\begin{eqnarray*}
(P_i^j)^TD_{jj}^kP_i^j&=&(P_i^j)^TA_{jj}^kP_i^j(P_i^j)^{-1}J_{j}(\lambda_i)P_i^j\\
&=&(P_i^j)^TA_{jj}^kP_i^jJ_{j}(\lambda_i)(P_i^j)^{-1}P_i^j\\
&=&\epsilon_i^j E_i^jJ_{j}(\lambda_i),
\end{eqnarray*}
where the second last equality is due to that $P_i^j$  commutes with $J_{j}(\lambda_i)$, see Lemma 2 in \cite{uhlig1976canonical}.
Thus the $j$th block of ($A_i, D_i$) is in the canonical form.
For the last $l-k$ elements in $\tilde{D }_i$, we have $$\tilde{A}_{(k+1):l,(k+1):l}={\rm diag}(a_{k+1},\ldots,a_l).$$ Let
$$P_i^{k+1}={\rm diag}(\frac{1}{\sqrt{|a_{k+1}|}},\ldots,\frac{1}{\sqrt{|a_l|}}),$$ then
$$(P_i^{k+1})^T\tilde{A}_{(k+1):l,(k+1):l}P_i^{k+1}={\rm diag}(\epsilon_i^{k+1},\ldots,\epsilon_i^l).$$ Denote
$P_i={\rm diag}(P_i^1,\ldots,P_i^k,P_i^{k+1})$, then $$\check{A}_i=P_i^T\tilde{A}_iP_i={\rm diag}(\epsilon_i^1 E_i^1,\ldots,\epsilon_i^k E_i^k,\epsilon_i^{k+1},\ldots,\epsilon_i^l)$$ is already in
the canonical form (\ref{CF1}). Besides,
 $$\check{D}_i=P_i^T\tilde{D}_iP_i=P_i^T\tilde{A}_i
C(\lambda_i)P_i=P_i^T\tilde{A}_iP_iC(\lambda_i)=\check{A}_iC(\lambda_i)$$ satisfies (\ref{CF2}) (the second last equality is due to that the two block diagonal matrices $C(\lambda_i)$ and $P_i$ commute).

Note in the above procedure, if there are too many $2\times 2$ blocks in the same Jordan chain, the time cost is very large. But this is really rare in practice since it is actually rare that any two eigenvalues of $A^{-1}D$ are the same and much rarer if two $2\times2$ blocks have the same eigenvalue, in which case  at least $4$ eigenvalues are the same.

For other Jordan chains, we only have $1\times 1$ Jordan blocks and thus it can be easily transformed to the canonical form by Algorithm 2 on
Lines 13--14 and 18--19 and the time required is very little when compared with the Jordan chains that has $2\times2$ blocks. \end{appendices}

\end{document}